

\magnification=\magstep1
\overfullrule=0pt
\font\tenmsb=msbm10
\font\sevenmsb=msbm7
\font\fivemsb=msbm5
\newfam\msbfam
\textfont\msbfam=\tenmsb
\scriptfont\msbfam=\sevenmsb
\scriptscriptfont\msbfam=\fivemsb
\def\Bbb#1{{\fam\msbfam\relax#1}}
\font\teneufm=eufm10
\font\seveneufm=eufm7
\font\fiveeufm=eufm5
\newfam\eufmfam
\textfont\eufmfam=\teneufm
\scriptfont\eufmfam=\seveneufm
\scriptscriptfont\eufmfam=\fiveeufm


\def\P{{\Bbb P}}     
\def\C{{\Bbb C}}     
\def\R{{\Bbb R}}     
\def\Q{{\Bbb Q}}     
\def\Z{{\Bbb Z}}     
\def\N{{\Bbb N}}     

\def\LL{{\cal L}} 
 
\def\AA{{\cal A}} 
\def\BB{{\cal B}} 
\def\CC{{\cal C}} 
\def\EE{{\cal E}} 
\def\HH{{\cal H}}
\def\II{{\cal I}}
\def\JJ{{\cal J}} 
\def\RR{{\cal R}}
\def\UU{{\cal U}}
\def\VV{{\cal V}}
\def\DD{{\cal D}}

\def\m2{m+2}
\def\noi{\noindent}

\def\linebreak{\relax\ifhmode\unskip\break\else
    \nonhmodeerr@\linebreak\fi}
\def\newline{\relax\ifhmode\null\hfil\break
\else\nonhmodeerr@\newline\fi}

\def\frac#1#2{{#1\over#2}}

\def\del#1{\frac {\partial} {\partial #1}}

\def\aa{\alpha}
\def\bb{\beta}

\centerline{\bf The $\AA$-hypergeometric
System}
\medskip
\centerline{\bf Associated with a Monomial Curve}
\footnote{}{\noi 1991 AMS Subject Classification: Primary 33C70, 
Secondary 14M05, 33D20.}
\medskip
\bigskip
\centerline{\bf Eduardo Cattani, Carlos D'Andrea, and Alicia Dickenstein}
\bigskip
\bigskip
\bigskip
\beginsection {\bf Introduction}

\bigskip

In this paper we make a detailed analysis of the
$\AA$-hypergeometric system (or GKZ system)
associated with a monomial curve and integral, hence resonant,
exponents. We describe all rational solutions and
show in Theorem~1.10 that they are, in fact,
Laurent polynomials. We also show that for any exponent,
there are at most two linearly independent 
Laurent solutions, and that the upper bound is reached if and
only if the curve is not arithmetically Cohen--Macaulay.
We then construct, for all integral parameters,
a basis of local solutions in terms
of the roots of the generic univariate polynomial (0.5) associated with
$\AA$. We also determine in Theorem~3.7 the holonomic rank $r(\aa)$ for all
$\aa\in\Z^2$ and  show that $d\leq r(\aa) \leq d+1$,
where $d$ is the degree of the curve. Moreover, the value $d+1$
is attained only for those exponents $\aa$ for which there
are two linearly independent rational solutions
and, therefore, $r(\aa) =d$ for all $\aa$ if and only if the
curve is arithmetically Cohen-Macaulay.

\smallskip

In order to place these results in their appropriate context,
we recall the definition of the $\AA$-hypergeometric systems.
These were introduced in
a series of papers in the mid  1980's by the Gel'fand
school, particularly Gel'fand, Kapranov, and
Zelevinsky  (see [7, 9], and the references
therein). Let $ \AA=\{\nu_1,\ldots,\nu_r\}\subset \Z^{n+1}$ 
be a finite subset which spans the lattice
$\Z^{n+1}$. Suppose, moreover, that there exists a vector
$\lambda = (\lambda_0,\ldots, \lambda_n) \in \Q^{n+1}$ such that
$\langle \lambda, \nu_j\rangle =1$ for all $j=1,\ldots,r$, \  i.e.
the set $\AA$ lies in a rational hyperplane.
Let  $\AA$ also denote the $(n+1)\times r$ matrix
whose columns are the vectors $\nu_j$.  Let $\LL\subset \Z^r$
be the sublattice of
elements $v \in \Z^r$ such that $\AA\cdot v = 0$.  Given 
$\aa\in \C^{n+1}$, the
$\AA$-{\sl hypergeometric system  with exponent (or parameter)} $\alpha$
is:
$$ \DD_v \varphi \ =\  0 \quad; \quad v\in \LL\eqno{(0.1)}$$
$$ \sum_{j=1}^r \nu_{ji} \, x_j 
\frac {\partial \varphi} {\partial x_j} \ =\   \alpha_i\,\varphi\quad ;
\quad i=1,\ldots,n+1\eqno{(0.2)}$$
where $\DD_v$ is
the differential operator in $\C^r$:
$$\DD_v \ :=\  \prod_{v_j > 0} \left(\del {x_j}\right)^{v_j} - 
\prod_{v_k < 0} \left(\del {x_k}\right)^{-v_k}\,.$$

The $\AA$-hypergeometric system is holonomic (with regular singularities) 
and,
consequently, the number of linearly independent solutions at
a generic point is finite [7].  Let $r(\aa)$  denote the
holonomic rank of the system, i.e. the
dimension of the space of local solutions at a generic point in $\C^r$.
If we drop the assumption that $\AA$ lies 
in a hyperplane, then the regular singularities property is lost
but, as Adolphson [2] has shown, the system remains holonomic.
The singular locus is described
by the zeroes of the principal $\AA$-determinant ([11]). 
We set 
$R := \C[\xi_1,\ldots,\xi_r]/\II_\AA$, where $\II_\AA$ is the
toric ideal 
$$\II_\AA\  := \ \bigl\langle\, \xi^u - \xi^v \ :\ u,v\in \N^r\,;
\ \AA\cdot u =
\AA \cdot v\,\bigl\rangle\,. \eqno(0.3)$$  
When $n=1$, we can assume without loss of generality that 
$$  \AA = \pmatrix{ 1&1&\cdots&1&1\cr 0&k_1&\cdots&k_m&d\cr } , 
\eqno (0.4) $$
where $0<k_1<\cdots k_m < d$. 
Note that the condition that the columns of $\AA$ generate the
lattice $\Z^2$ is equivalent to $\gcd(k_1,\ldots,k_m,d)=1$.
The homogeneous ideal $\II_\AA$ defines 
a monomial
curve $X_\AA \subset \P^{m+1}$ of degree $d$ whose homogeneous coordinate ring
is $R$. $X_\AA$ is {\sl normal} if and only if $d=m+1$.  Recall that
$X_\AA$ is said to be arithmetically Cohen-Macaulay if and only if
the ring $R$ is Cohen-Macaulay. 

The system associated with (0.4) admits very interesting solutions. 
Let 
$$ f(x;t) := 
x_0 + x_{k_1}\,t^{k_1} + \cdots +  x_{k_m}\,t^{k_m} + x_d\,t^d\ , 
\eqno(0.5)$$
denote the generic polynomial with exponents $0,k_1,\ldots,k_m,d$.
It is not hard to see that
the powers  $\rho^s(x)$, $s\in\Z$, of the roots of $f(x;t)$,
viewed as functions of the coefficients,
are algebraic solutions of the $\AA$-hypergeometric system with exponent 
$(0,-s)$. This fact was observed by Mayr [17] who
constructed series expansions for these functions.  These have
more recently been refined by Sturmfels [22].
The total sum
$$p_s(x)\ :=\ \rho_1^s(x) + \cdots + \rho_d^s(x)\eqno(0.6)$$
will then be a rational solution with the same exponent.  Similarly,
one can show that the local residues
$${\rm Res}_{\rho(x)}\,\left( \frac {t^b} {f^a(x;t)}
\frac {dt} {t} \right)\ ;\quad
a,b\in\Z, \ a\geq 1\eqno{(0.7)}$$
give algebraic solutions with exponent $(-a,-b)$ and, again,
the total sum of residues will be a rational solution.

In \S1 we describe explicitly all rational solutions
of the $\AA$-hypergeometric system associated with
a monomial curve.  Since for $\AA$ as in (0.4), the principal 
$\AA$-determinant  factors into powers of $x_0$,
$x_d$, and the discriminant $\Delta(f)$, we know
a priori what the possible denominators of a rational
solution may be.  
However, we show in Theorem~1.10 that there are no rational solutions
whose denominator involves $\Delta(f)$ and therefore,
every rational solution must be a
Laurent polynomial.  This is a somewhat surprising result
which is peculiar to the case $n=1$ (see Example 1.11).
One may give explicit
formulas for these Laurent polynomials in terms of hypergeometric
polynomials in fewer variables.  When applied to the sum of powers
of roots one recovers the classical Girard formulas. 
One also obtains similar expressions for 
total residues in terms of hypergeometric polynomials.

We show that for any  $\alpha\in\Z^2$ 
the dimension of the space of rational $\AA$-hypergeometric
functions with parameter $\alpha$ is at most $2$.  Moreover,
the value $2$ may be reached for only finitely
many values of $\aa$ and this happens  if and only if the
ring $R$ is not Cohen-Macaulay.  

In \S2 we exhibit a family of algebraic $\AA$-hypergeometric
functions defined  in terms of the roots of the polynomial $f(x;t)$.
These are the building blocks for the construction,
in \S 3, of local bases of solutions  and the determination of 
the holonomic rank for all integral exponents.
It becomes necessary to consider four possibilities
for the exponent $\aa$.  These cases admit combinatorial
descriptions (see (1.9)) and correspond to the existence of a polynomial
solution; a one-dimensional space of rational -non-polynomial- solutions; 
a two-dimensional space of rational solutions;
or no rational solution for the given exponent.
For $u\in\N^{m+2}$, the derivative $D_u$ maps 
$\AA$-hypergeometric functions to $\AA$-hypergeometric
functions while changing the exponent from $\aa$ to
$\aa - \AA\cdot u$.
A careful analysis of the kernel and image of this operator
together with Corollary 5.20 of [2] leads to the determination
of the holonomic rank for all values of $\aa$.
We show, in particular, that
$$d \leq r(\aa) \leq d + 1 \eqno(0.8)$$
and that $r(\aa) = d+1$ exactly for
those parameters $\alpha\in\Z^2$ for which the dimension 
of the space of rational solutions is $2$.  
Hence, $r(\aa) = d$ for all $\aa\in\Z^2$ if
and only if the curve $X_\AA$ is arithmetically
Cohen-Macaulay.

These results allow
us to clarify the relationship between 
the holonomic rank and ${\rm vol}(P)$, the
normalized volume  of the convex hull
$P\subset \R^{n+1}$ of $\AA$ and
the origin (which equals the degree $d$ in the case of curves).
It was originally claimed in [9,~Theorem~2]
that $r(\aa) = {\rm vol}(P)$ in all cases, but it
was pointed out by Adolphson that, for
resonant exponents, 
the proof required the assumption that the ring $R$ 
be Cohen-Macaulay (see [10]). 
In [2,~Corollary~5.20], Adolphson showed that 
$r(\aa) = {\rm vol}(P)$, for $\aa$ a semi-nonresonant exponent,
without any additional assumptions on $R$.
The first explicit example where the equality fails is given
in [23] and described in Example~1.8.i).  
In the forthcoming monograph [19], Saito, Sturmfels and Takayama
prove (0.8) using Gr\"obner deformation methods and show that
the inequality 
$r(\aa) \geq {\rm vol}(P)$ holds without restrictions on $n$.

\medskip

Very little seems to be known about the problem of finding
rational solutions of differential equations beyond the case
of linear differential operators in one variable.  In this case,
Singer [20] has shown that one can  determine 
in a finite number of steps whether a
given equation has a rational solution and find a basis for
the space of such solutions.  Abramov 
and Kvasenko ([1]) have further studied 
the problem of effectively finding rational solutions for
such operators.  In our case one
could, in principle, use non-commutative elimination to
obtain linear operators in one variable with coefficients
that depend rationally on the other variables and apply
Singer's decision procedure to characterize the rational
solutions.  
This can be done in small examples using 
non-commutative Gr\"obner bases packages such as {\tt kan} ([24])
but we have not been able to obtain any general results 
in this manner.

Gel'fand, Zelevinsky,
and Kapranov have constructed series solutions for (0.1)-(0.2), 
associated with regular triangulations
of the polytope $P$.  When the exponent $\alpha$ is non-resonant
for the triangulation, it is possible to obtain in this manner
${\rm vol}(P)$-many independent solutions.  There are, however,
very interesting cases in which the exponents are integral and,
therefore, automatically resonant.  For example, it has been
observed by Batyrev in [4] that
the period integrals of Calabi-Yau hypersurfaces in toric varieties
satisfy an $\AA$-hypergeometric system with exponents 
$\alpha=(-1,0,\ldots,0)$.  In this case, series solutions have
been obtained by Hosono, Lian, and Yau in [14, 15] 
(see also  [3, 5, 13]).  
Very recently, Stienstra [21] has generalized the $\Gamma$-series
construction of Gel'fand, Zelevinsky,
and Kapranov to obtain series solutions in the case of
resonant exponents under a maximal-degeneracy
assumption.  In particular, if $\aa=0$ and $P$
admits a unimodular triangulation --which implies that $R$ is
Cohen-Macaulay-- all solutions of (0.1)-(0.2) may be obtained 
in this manner.  This method also yields all solutions
of interest in the context of toric mirror symmetry.

A key result, in the curve case, is Theorem~1.9 which asserts 
that for $m\geq 1$, there are no rational
solutions with integral exponents in the Euler-Jacobi cone.
 This corresponds to the
classical vanishing theorem for the total sum of residues, a
statement which has a generalization as the
Euler-Jacobi theorem (see [16] for the most general
form of this result). It is interesting to note then that
Euler-Jacobi vanishing is  a  consequence of
the fact that residues satisfy the $\AA$-hypergeometric 
system. We also point out that while the characterization of
Laurent solutions follows from formal arguments, the proof of
the Euler-Jacobi vanishing involves  transcendental
methods.

\medskip\noi
{\bf Acknowledgments:} 
We are grateful to 
B.~Sturmfels for many helpful comments and, particularly, 
for the statement and proof of Proposition~1.6.
We also thank A.~Zelevinsky for useful suggestions
and the two referees for their careful reading of the manuscript
and very thoughtful suggestions for improvement.
E.~Cattani was supported by NSF Grant DMS-9406462.
C.~D'Andrea was supported by a Doctoral Fellowship
from FOMEC, Argentina.  A.~Dickenstein was supported by UBACYT
and CONICET, Argentina.

\bigskip
\bigskip

\beginsection {\bf 1. Rational solutions}

\bigskip

The polynomial solutions of a
general
$\AA$-hypergeometric system 
admit a very simple description.  
Given $\alpha\in \Z^r $, we define  the
{\it hypergeometric polynomial}
$$ \Phi^{\AA}( \alpha; x) \  := \  
\sum_
{ \scriptstyle u\in \N^r \atop \scriptstyle \AA\cdot u=\alpha} 
\ {x^u \over u \, !} \  = \  
\sum_{ \scriptstyle u\in \N^r \atop \scriptstyle 
\AA\cdot u=\alpha }
\ { x_1^{u_1} x_2^{u_2} \cdots x_r^{u_r} \over
u_1 ! \, u_2 ! \cdots u_r ! }. \eqno (1.1)$$
As usual, we set 
$\Phi^{\AA}( \alpha; x) := 0$ if $\alpha\not\in \AA\cdot \N^r$.
The following result, whose verification is left to the reader, is
Proposition~2.1 in [18]:
 
\medskip
\proclaim Proposition 1.1. $\Phi^{\AA}( \alpha; x)$ is the unique,
up to scaling, polynomial solution of the 
$\AA$-hypergeometric system with exponent $\alpha$.  Moreover, for
any $u\in \N^r$,
$$ D_u \bigl( \Phi^{\AA}( \alpha; x) \bigr) = 
\Phi^{\AA}( \alpha - \AA\cdot u; x) \eqno(1.2)$$
where $D_u$ stands for the partial derivative 
$\partial^{|u|}/\partial x^u$.

\medskip

The purpose of this section is to describe the rational solutions of the 
$\AA$-hyper\-geometric 
system associated with a matrix $\AA$ as in (0.4).  Note that
for $m=0$ the system restricts to the homogeneity equations
(0.2). Therefore we may assume throughout that $m\geq 1$. 
To simplify our notation, we will index all $(m+2)$-tuples by 
$0,k_1,\ldots,k_m,d$.  Let $e_0,e_{k_1},\ldots,e_d$ denote
the standard basis of $\Z^{m+2}$.  For $i=1,\ldots,m$ we have
$$\omega_{k_i}\  :=\  (d-k_i)\,e_0 - d\, e_{k_i} + k_i\, e_d\in \LL\,.
\eqno {(1.3)}$$

The following 
observation will be  useful in the sequel:

\medskip
\proclaim Proposition 1.2.  Suppose $\varphi$ is a local
holomorphic
solution of
(0.1)-(0.2), polynomial with respect to any of the variables
$x_0,x_{k_1},\ldots,x_d$.  Then $\varphi$ is a Laurent polynomial. 

\smallskip\noi
{\bf Proof:} Since $\varphi$ satisfies the equations (0.1), 
it follows from (1.3) that
 for all $\ell\in\N$, 
$$D_{k_i}^{\ell d} \varphi  \ =\ D_{0}^{\ell (d-k_i)} 
D_{d}^{\ell k_i}\varphi\ 
\eqno{(1.4)}$$
and, consequently, if $\varphi$ is polynomial in any of the
variables, it must be so in all of the variables
$x_{k_i}$, $i=1,\ldots,m$,  and we may write:
$$\varphi(x) \ =\ \sum_u\,\varphi_u(x_0,x_d)\,\check x^u
\quad;\quad \check x =(x_{k_1},\ldots,x_{k_m})\,
$$
where $u$ varies in a finite subset of $\N^m$ and each
$\varphi_u(x_0,x_d)$ is homogeneous in each variable with
respective degrees $\bb_0$, $\bb_d$ satisfying
$d\beta_d \in \Z$, $\bb_0 + \bb_d\in \Z$.  Thus,
$\varphi_u(x_0,x_d) = c_u\,x_0^{\bb_0}\,x_d^{\bb_d}$, $c_u\in \C$.
But, because of (1.4), for $r,s\in\N$ sufficiently large,
$D_{0}^{r} 
D_{d}^{s}\varphi =0$,
and consequently one, and therefore both, of $\bb_0, \bb_d$ must
be an integer.  Hence $\varphi$ is a Laurent polynomial.\ \ $\diamond$

Note that in the one-dimensional case, an $\AA$-hypergeometric
Laurent polynomial may not contain a 
non-zero term of
the form $c_u\,x^u$ with $u_{k_i}<0$.
This follows from the fact that the singular locus $\Sigma$ of
the hypergeometric system is given by the zeroes
of the principal $\AA$-determinant, i.e.
$$\Sigma \ =\  \{ x_0 =0\}\ \cup\  \{ x_d =0\}\ \cup\ 
\{ \Delta(f) =0\}$$
where $\Delta(f)$ is the discriminant of the generic polynomial (0.5).
Alternatively, if a solution contains 
a non-zero term $c_u\,x^u$ with $u_{k_i}<0$,
being
in the kernel of the differential operator $D_{{k_i}}^d - D_0^{d-k_i}
D_d^{k_i}$,
it must also contain non-trivial terms of the form 
$c_v\,x^v$ with $v_{k_i} = u_{k_i} - j\,d$ for all positive integers
$j$ and this is clearly
impossible.   A similar argument shows that 
a Laurent solution may not contain terms of the form:
$ c_u\, {x^u} $
with both $u_0<0$ and $u_d<0$.

Thus, any Laurent solution must be of the form
$$L_0(x) + L_d(x)$$
where $L_0(x)$ has as denominators only powers of $x_0$ and 
$L_d(x)$ has as denominators only powers of $x_d$.  We note that 
the study of $L_0(x) $ and $ L_d(x)$ is completely symmetric.
Indeed, let $\ell_j = d-k_{m-j+1}$, $j=1,\ldots,m$, and
$$  \hat \AA = \pmatrix{  1&1&\cdots&1&1\cr 0&\ell_1&\cdots&\ell_m&d\cr
}\eqno{(1.5)}$$ 
Then, if $R(x)$ is a 
Laurent solution of the $\AA$-hypergeometric system and exponent
$\alpha=(\aa_1,\aa_2)$
(although $\aa$ should be viewed as a column vector
we will, for simplicity of notation, always write exponents
as row vectors), the function $\hat R(y_0,y_{\ell_1},\ldots,y_{\ell_m},y_d)$
obtained from $R$ by substituting: 
$$x_0 \mapsto y_d\ , \quad x_{k_j} \mapsto y_{\ell_j}
\ , \quad x_d \mapsto y_0, \eqno{(1.6)}$$
is a solution of the $\hat \AA$-hypergeometric system and exponents
$\hat\alpha=(\alpha_1, d\alpha_1-\alpha_2)$.

\bigskip

\proclaim Lemma 1.3. 
For $\alpha \in \AA\cdot\N^{m+2}$, the only Laurent solutions of the 
$\AA$-hypergeometric system are the constant multiples of the
hypergeometric polynomial (1.1).
 
\noi{\bf Proof:} Suppose there is a Laurent solution $L(x)$ of exponent 
$\alpha$ containing a term of the form $x^u/x_d^r$, with $u\in\N^{m+2}$,
$u_d=0$ and
$r>0$ (we will always assume that monomials are written in
reduced form). Then $\AA\cdot (u - r e_d) = \alpha$. 
Let $v\in  \N^{m+2}$ be such that $\AA\cdot v = \alpha$, then
$\AA\cdot u =\AA\cdot (v + r\,e_d) $,  and
the operator $D_u - D_d^rD_v$
being  in the hypergeometric system, must vanish on $L$.
This means that $L$ must also contain a term $x^w$ whose
derivative $D_d^rD_v(x^w)$
is a non-zero multiple of $\,1/x_d^r$.  Since $v_d \geq 0$, this
is clearly impossible.  Arguing by symmetry, we see that there cannot
be a solution containing a term of the form 
$\,x^u/x_0^r$, with $u_0=0$ and
$r>0$. $\ \ \diamond$
 
\bigskip
 
We will denote by $\BB$  the matrix
$$  \BB = \pmatrix{  1&1&\cdots&1\cr 0&k_1&\cdots&k_m\cr}$$ 
and  $x'$ the vector consisting of the first  
$m+1$ variables $(x_0,x_{k_1},\ldots,x_{k_m})$.  
Similarly,  let $\CC$ be the matrix
$$  \CC = \pmatrix{  1&\cdots&1&1\cr k_1&\cdots&k_m&d\cr}$$ 
and $\tilde x$ the vector $(x_{k_1},\ldots,x_{k_m},x_d)$.

\bigskip

Given $\alpha\in\Z^2$ we define
$$\Psi^{\AA}_d(\alpha;x) 
\ := \ \sum_{r\geq 1} \ (-1)^r \, (r-1)!\  \frac { \Phi^{\BB}( \alpha'(r); 
x')} {x_d^r}\ , \eqno {(1.7)}$$
where $\alpha'(r) = \alpha + r\,(1,d)$, and
$$ \Psi^{\AA}_0(\alpha;x) 
\ := \  \sum_{r\geq 1} \ (-1)^r \, (r-1)!\  \frac { \Phi^{\CC}
( \tilde\alpha(r); 
\tilde x)} {x_0^r}\ , \eqno {(1.8)}$$
where $\tilde\alpha(r) = \alpha + r\,(1,0)$.  

Note that both sums are finite.  
This follows from the fact that 
${ \Phi^{{\BB}}( \alpha'(r); 
x')} =0$ unless $\alpha'(r) \in {\BB}\cdot\N^{m+1}$, but then
 $\alpha_2+dr\leq k_m(\alpha_1+r)$ and, therefore, 
$(d-k_m)r \leq k_m\alpha_1- \alpha_2$.  
This means that the possible values of $r$ in (1.7)
are bounded by $(k_m\alpha_1-\alpha_2)/(d-k_m)$.   The statement for
(1.8) follows by symmetry.  In fact, we should observe
that the change of variables $(1.6)$ transforms
$\Psi_d^\AA(\alpha;x)$ into $\Psi_0^{\hat\AA}(\hat\alpha;x)$.

The following subsets of $\Z^2$ will play an important role in
the description of $\AA$-hyper\-geometric functions:
$$
I(\AA)\ := \ \AA\cdot \N^{m+2}\ ;$$
$$F_0(\AA) \ :=\  \CC\cdot \N^{m+1} - \N\, (1,0) \ ;
\quad E_0(\AA) \ :=\  F_0(\AA) \setminus I(\AA)\ ;$$
$$F_d(\AA) \ :=\  \BB\cdot \N^{m+1} - \N\,  (1,d) \ ;
\quad E_d(\AA) \ :=\  F_d(\AA) \setminus I(\AA)\ ;\eqno(1.9) $$
$$
E(\AA) \ := E_0(\AA) \cap E_d(\AA)\ ;\quad
J(\AA) \ :=\ \Z^2 \setminus (I(\AA)\cup F_0(\AA) \cup F_d(\AA)).$$
Note that via the change of variables (1.6), and denoting 
for $\aa=(\aa_1,\aa_2)\in\Z^2$,
$\hat\aa = (\aa_1,d\aa_1-\aa_2)$, we have for $\hat\AA$ as in (1.5):
$$I(\hat\AA) = \widehat {I(\AA)}\ ;\quad
F_0(\hat\AA) = \widehat{F_d(\AA)}\ ;\quad
F_d(\hat\AA) = \widehat{F_0(\AA)}\,. $$
It is clear from the definitions (1.7) and (1.8) that 
$\Psi^{\AA}_0(\alpha;x)\not=0$ if and only if $\aa\in F_0(\AA)$
and, similarly, 
$\Psi^{\AA}_d(\alpha;x)\not=0$ if and only if $\aa\in F_d(\AA)$.
In particular,
$$\Psi^{\AA}_d(\alpha;x)=0\quad \hbox{if }
\quad d\,\alpha_1 < \alpha_2
\,,\eqno(1.10)$$
and, 
$$\Psi^{\AA}_0(\alpha;x)=0\quad \hbox{if }
\quad  \alpha_2 < 0
\,.\eqno(1.11)$$
On the other hand, the importance of the sets $E_0(\AA)$ and
$E_d(\AA)$ stems from the fact that according to Lemma 1.3, it
is only for $\aa\in E_0(\AA)$ (respectively $\aa\in E_d(\AA)$)
that the Laurent polynomial $\Psi^{\AA}_0(\alpha;x)$
(respectively $\Psi^{\AA}_d(\alpha;x)$) may be -and as the
following result shows is- $\AA$-hypergeometric.  
Note also that there are no $\AA$-hypergeometric Laurent polynomials
with exponent $\aa\in J(\AA)$ and it will be a consequence of 
Theorem 1.10 that there are no rational 
$\AA$-hypergeometric functions with exponent $\aa\in J(\AA)$.

\bigskip

\proclaim Theorem 1.4. Let $\AA$ be as in (0.4) and 
$\alpha \in \Z^2 \setminus I(\AA)$.
Then:
\smallskip\noi
{i)} $\Psi^{\AA}_d(\alpha;x)$ and $\Psi^{\AA}_0(\alpha;x)$
are solutions (possibly trivial)
of the $\AA$-hypergeometric system with parameter 
$\alpha$.
\smallskip
\noi
{ii)} For any $u\in\N^{m+2}$,
$$ D_u \bigl( \Psi^{\AA}_d( \alpha; x) \bigr) \ = \ 
\Psi^{\AA}_d( \alpha - \AA\cdot u; x) \ ;
\quad 
D_u \bigl( \Psi^{\AA}_0( \alpha; x) \bigr) \ = \   
\Psi^{\AA}_0( \alpha - \AA\cdot u; x)\eqno(1.12)$$
\smallskip
\noi
{iii)} The functions $\Psi^{\AA}_d(\alpha;x) $ and 
$\Psi^{\AA}_0(\alpha;x) $ span the space of 
Laurent solutions of the $\AA$-hypergeometric system with parameter
$\alpha$.
 
\bigskip
  
\noi{\bf Proof:} \ Clearly, i) is an immediate consequence of ii) and,
because of symmetry, it suffices to show (1.12) for   
$\Psi^{\AA}_d(\alpha;x)$.
Suppose $u\in \N^{m+2}$ is such that $u_d=0$, then
$$ \eqalign{
D_u \bigl( \Psi^{\AA}_d( \alpha; x) \bigr) \quad &= \quad
\sum_{r\geq 1} \ (-1)^r \, (r-1)!\  \frac 
{ D_u\bigl(\Phi^{{\BB}}( \alpha'(r); x')\bigr)} {x_d^r}\cr 
&=\quad 
\sum_{r\geq 1} \ (-1)^r \, (r-1)!\  \frac 
{ \Phi^{{\BB}}( \alpha'(r) - {\BB}\cdot u; 
x')} {x_d^r} \quad = \quad \Psi^{\AA}_d( \alpha - \AA\cdot u; x)
\,.\cr}$$
Thus, it remains to prove (1.12) for the partial derivative
$D_d$.  We have
$$\eqalign{
D_d\,\Psi^{\AA}_d( \alpha; x) \quad
&= \quad
\sum_{r\geq 1} \ (-1)^{r+1} \, r!\  \frac 
{ \Phi^{{\BB}}( \alpha'(r); x')} {x_d^{r+1}}\cr
&= \quad
\sum_{r\geq 1} \ (-1)^{r+1} \, r!\  \frac 
{ \Phi^{\BB}( (\alpha-\AA\cdot e_d)'(r+1); x')} {x_d^{r+1}}\cr
&=\quad
\Psi^{\AA}_d( \alpha - \AA\cdot e_d; x)\,.\cr}
$$ 
The last equality follows since $\alpha\not\in \AA\cdot 
\N^{m+2}$ implies that 
$$ \Phi^{\BB}( (\alpha-\AA\cdot e_d)'(1); x')\quad=\quad
\Phi^{\BB}( \alpha; x')\quad=\quad0\,.$$  

\medskip
 
\noi Suppose now that $L(x)$ is a Laurent solution with exponent
$\alpha$ and write $L(x) = L_d(x) + L_0(x)$.  If we decompose
further:
$$L_d(x) = \sum_{r\geq 1} \frac {A_r(x')} {x_d^r}$$
then, the polynomials $A_r(x')$ must be solutions of the
$\BB$-hypergeometric system and exponent $\alpha'(r)$. 
Thus, $A_r(x') = c_r\, \Phi^{\BB}(\alpha'(r);x')$.  Assume now
that for $r,s\geq 1$, $r\not= s$, we have $\aa'(r),\aa'(s)\in I(\BB)$
and
let  $v,w\in \N^{m+1}$ be 
such that $\BB\cdot v =\alpha'(r)$ and $\BB\cdot w=\alpha'(s)$.  
Then, since
$\AA\cdot (v + se_d) = \AA\cdot (w + re_d)$, the operator
$$\frac {\partial^{|v|+s}} {\partial x_{}^v\,x_d^s}\  -\ 
\frac {\partial^{|w|+r}} {\partial x_{}^{w}\,x_d^{r}}$$
is in the hypergeometric system and must vanish on $L$.  This means
that 
$$\frac {(-1)^s\, (r+s-1)!} {(r-1)!\,{x_d^{r+s}}} 
\ \frac {\partial^{|v|} (A_r(x')) } {\partial x_{}^v}\ 
 =\ 
\frac {(-1)^r\, (r+s-1)!} {(s-1)!\,{x_d^{r+s}}}
\ \frac {\partial^{|w|} (A_s(x')) } {\partial x_{}^w}\ 
$$
which implies that
$$\frac {(-1)^s } {(r-1)!} 
\ c_r
 =\ 
\frac {(-1)^r} {(s-1)!}  
\ c_s
$$
and, consequently, 
$$c_r = c \ (-1)^r \ (r-1)!\,.$$
A symmetric argument shows that if the component $L_0(x)$ 
is non-trivial
then it must be a constant multiple of 
$\Psi^{\AA}_0(\alpha;x)$, which
proves part iii).
\ \ $\diamond$

\medskip

We state for emphasis: 
\bigskip
\proclaim Corollary 1.5. Let $\AA$ be as in (0.4) and 
let $\LL(\aa)$ denote the vector
space of $\AA$-hyper\-geometric Laurent polynomials of exponent
$\aa$.  Then $\dim \LL(\aa) \leq 2$.  Moreover,
$\dim \LL(\aa) \geq 1$  if and only if 
$\aa\in I(\AA)\cup E_0(\AA) \cup E_d(\AA)$,     and 
$\dim \LL(\aa) = 2$ if and only if $\aa\in E(\AA)$.

\medskip
The following result was brought to our attention
by Bernd Sturmfels:
 
\bigskip

\proclaim Proposition~1.6. Given $\AA$ as in (0.4), 
the ring $\ R = \C[\xi_0,\ldots,\xi_d]/\II_\AA\ $ is 
Cohen-Macaulay if and only if $E(\aa) =\emptyset$.

\medskip\noi
{\bf Proof:\ } The ring $R$ is a particularly simple
example of an affine semigroup ring whose properties
have been extensively studied (see for example
[6,~Chapter 6], [12], [25]).
In fact, Proposition~1.6
is a special case of Theorem~2.6 in [12]
which gives necessary and sufficient conditions for an
affine semigroup ring which, like $R$ does, admits a
system of monomial parameters.  Their condition (ii)
is easily seen to be equivalent, in our
notation, to $E(\aa) =\emptyset$.   \ \ $\diamond$  

\bigskip\noi
{\bf 1.7. Remarks: } i) When the curve $X_\AA$ is normal,
i.e. $d=m+1$, the sets defined in (1.9) have a very simple
description:  the image $I(\AA)$ coincides with the ``cone"
(properly speaking semigroup)
$$\CC := \{\aa=(\alpha_1,\alpha_2)\in\Z^2\,:\,
0\leq \alpha_2 \leq d\alpha_1\}\,,$$
while
$E_0(\AA) = \{\aa\in\Z^2 : 
\alpha_2 \geq 0, d\alpha_1 < \aa_2\}$, and
$E_d(\AA) = \{\aa\in\Z^2 : 
\alpha_2 < 0, d\alpha_1 \geq \aa_2\}$.  
The complement $J(\AA)$ of these
three sets, i.e. the set of $\aa\in\Z^2$ for which there are
no $\AA$-hypergeometric Laurent polynomials of exponent 
$\aa$ is the {\sl Euler-Jacobi cone}:
$$\EE\JJ := \{(\alpha_1,\alpha_2)\in\Z^2\,:\,
d\alpha_1<\aa_2<0\}\,.\eqno{(1.13)}$$
\smallskip

\noi ii) For arbitrary $\AA$ as in (0.4) we have
$I(\AA)\subset \CC$, $E(\AA)\subset \CC$, and
$\EE\JJ\subset J(\AA)$.

\smallskip
\noi iii)\  The conditions (1.10) and (1.11) are far from being
sharp.  It is easy to see, for example, that for
$\aa\in\CC$, $\Psi_0^\AA(\aa;x)\not=0$ and 
$\Psi_d^\AA(\aa;x)\not=0$, imply that 
$\ k_1\,\aa_1 < \aa_2 < k_m\,\aa_1$.
In particular, if
$m=1$, $k_1=k_m$ and  there  is no $\alpha\in \CC$ for which  both 
$\Psi^{\AA}_0(\alpha;x)$ and $\Psi^{\AA}_d(\alpha;x)$ are non-trivial.
Hence, there is for each $\alpha$ at most one, up to constant
multiple, Laurent solution of exponent $\alpha$.  Note also that
the toric ring $R$ is always Cohen-Macaulay but it is normal if
and only if $d=2$.

\smallskip
\noi iv)\  It is not hard to prove (see for example [2,~Lemma~3.12]) 
that there exists $v\in \AA\cdot \N^{m+2}$ such that 
$v + \CC \subset I(\AA)$.  Thus, for $\aa \in v + \CC$, there
is a unique $\AA$-hypergeometric Laurent polynomial and it is given
by (1.1).  Moreover, when this observation is combined with the
inequalities in iii), it follows that the set $E(\AA)$ is finite.

\bigskip
\noi {\bf 1.8. Examples:} \ We exhibit two examples where $E(\AA)
\not= \emptyset$ and, consequently, the associated toric ring $R$ is not
Cohen-Macaulay.
\medskip
\noi i) This is the ``running example" in [23]. Let 
$$\AA = \pmatrix{1&1&1&1\cr
0&1&3&4\cr}$$
The exponent $\aa=(1,2)$ is the unique element in $\CC$ such that 
$\aa\not\in I(\AA)$ and $\aa_1 < \aa_2 < 3\aa_1$.  
Both
$\Psi^{\AA}_0((1,2);x) = (-1/2)\,(x_1^2/x_0)$ and
 $\Psi^{\AA}_4((1,2);x) = (-1/2)\,(x_3^2/x_4)$ are non-trivial.

\noi The element $\beta =(2,3)\in F_0(\AA) \cap F_4(\AA)$ and
therefore both $\Psi^{\AA}_0(\bb;x) =(-1/6)\,x_1^3/x_0$ and
$\Psi^{\AA}_4(\bb;x) =(-1/2)\,x_1x_3^2/x_4$ are non-trivial.
However, since $\beta = \AA\cdot(e_0 + e_3) \in I(\AA)$, it follows
from Lemma~1.3 that neither $\Psi^{\AA}_0(\bb;x)$,
nor $\Psi^{\AA}_d(\bb;x)$, is $\AA$-hypergeometric and the
only $\AA$-hypergeometric Laurent polynomials of exponent $\bb$
are the multiples of the polynomial
$\Phi^\AA(\bb;x) = x_0x_3$.

\medskip
\noi ii) Consider the system associated with the matrix
$$\AA = \pmatrix{1&1&1&1&1\cr
0&6&7&13&14\cr}$$
and $\aa=(2,18)\not\in I(\AA)$.  We have 
$$\Psi^{\AA}_0((2,18);x) \ = \ -\frac {1} 6 \, 
\frac {x_6^3} {x_0}\quad ;\quad
\Psi^{\AA}_{14}((2,18);x) \ = \ -\frac {1} 2 \, 
\frac {x_6x_{13}^2} {x_{14}}+
\frac 1 6 \, \frac {x_7x_{13}^3} {x_{14}^2}\ .$$

\bigskip

It follows from (1.10) and (1.11) that 
there are no $\AA$-hypergeometric Laurent polynomials whose
exponent $\alpha$ is in 
the Euler-Jacobi cone (1.13).
In fact, as the 
following result shows, there are no
{\sl rational} solutions with exponent in that region.
 
\bigskip

\proclaim Theorem 1.9. 
The $\AA$-hypergeometric system associated with the matrix (0.4) has
no rational solutions whose exponent $\alpha$ lies in the Euler-Jacobi cone.
 
 \bigskip
  
Given $\aa\in\Z^2$, we will denote by $\RR(\aa)$ the vector
space of rational $\AA$-hypergeometric functions of exponent
$\aa$.
Before proving Theorem~1.9 we note the following consequence.
 
\bigskip
\proclaim Theorem 1.10.
The only rational solutions of the 
$\AA$-hypergeometric system associated with the matrix (0.4)
are the Laurent polynomial solutions described by 
Proposition~1.1  and Theorem~1.4.
 
\medskip
\noi {\bf Proof:\ } Suppose $\varphi\in\RR(\aa)$.
For $\ell\in\N$ sufficiently large, 
$\beta = \aa - \ell \,(1,k_1)$ lies in the Euler-Jacobi cone.
Then,
 $D^\ell_{k_1}\varphi$ is a rational solution in  
$\HH(\beta)$, and thus it is identically zero by Theorem~1.9.
Hence, $\varphi$
is polynomial in $x_{k_1}$ and by Proposition~1.2 it must be a
Laurent polynomial.\ \ $\diamond$

\bigskip
\noi{\bf 1.11.\ Example: }  Theorem 1.10 is not true
for $n>1$.
It fails already in the simplest
two-dimensional situation: 
Consider   the hypergeometric
system associated with the matrix
$$\AA = \pmatrix{1&1&1&1\cr
0&1&0&1\cr
0&0&1&1\cr}\,.$$
The lattice $\LL$ has rank one; in fact $\LL = \Z\cdot (1,-1,-1,1)^T$, 
and the system (0.1)-(0.2) is equivalent
to Gauss' classical hypergeometric equation.  
The function $1/(x_1x_4-x_2x_3)$ is a solution  with
parameters $(-2,-1,-1)$.

\bigskip

\noi {\bf Proof of Theorem~1.9:\ }  The proof will be by
induction on $m$.  We begin by considering the case $m=1$
and write, for simplicity, $k_1=k$.  Note that in this
case the lattice $\LL$ has rank $1$ and is
generated by $\omega:=(d-k)e_0 - de_k + k e_d$.
In particular, for appropriate values of $\aa$, we can write the
$\AA$-hypergeometric functions in terms of classical hypergeometric
functions (see [9,~\S3.1]).
 
\smallskip
\noi The discriminant of the generic polynomial 
$x_0 + x_k\,t^k + x_d\,t^d$ is, up to factors which are powers
of $x_0$ and $x_d$,
$$\Delta(x) \ =\  
d^d\, x_0^{d-k} \,x_d^k + (-1)^{d-1}\, k^k\, (d-k)^{d-k}\, x_k^d 
 \ =\  {c\,x_k^d}\, {(1-\lambda\, z)}$$
where 
$\ c  \ =\ (-1)^{d-1} \,k^k\, (d-k)^{d-k}$, 
$\ z  =
x^\omega = x_0^{d-k}\,x_k^{-d}\,x_d^k $,
and
$$\lambda \ =\  \frac {(-1)^d\,d^d} {k^k\, (d-k)^{d-k}}$$
Suppose now that $R(x) = P(x)/Q(x) \in\RR(\aa)$ with
$\aa\in\EE\JJ$, i.e. $d\alpha_1<\alpha_2<0$.  Note that both $P$ and $Q$ are
bihomogeneous relative to the $\Z^2$-degree defined by $\AA$.
We can then write $P(x) = x^v\, P_1(z)$, where $P_1$ is a polynomial
and $P_1(0)\not=0$.
Thus, up to a constant,
$$R(x)  \ =\  x^u \ \frac {P_1(z)} {(1-\lambda\,z)^r}
 \ =\ x^u \,\sum_{j\geq 0}\  c_j\, z^j
 \ =\ x^u \,\sum_{j\geq 0}\  c_j\, x^{j\omega}\,,
\eqno {(1.14)}$$
with $ c_0\not=0$, and $\AA\cdot  u=\alpha$.  
\smallskip\noi
Since $R$ is
in the kernel of the differential operator 
$$\DD_\omega\ =\ D_0^{d-k}D_d^k - D_k^d\,, \eqno {(1.15)}$$
the coefficients in (1.14)  satisfy the relation:
$$
c_j\prod_{\ell=0}^{d-1} (u_k - jd-\ell)=
c_{j+1}\prod_{\ell'=0}^{d-k-1} (u_0 + (j+1) (d-k) -\ell')
 \prod_{\ell''=0}^{k-1} (u_d + (j+1) k -\ell'')\,,
\eqno(1.16)$$
for all $j\in\Z$. Setting $j=-1$ we get
$$0 = 
c_{0}\prod_{\ell'=0}^{d-k-1}(u_0 -\ell')
\prod_{\ell''=0}^{k-1} (u_d  -\ell'')\,,$$
which implies that either $0\leq u_0 \leq d-k-1$, or 
$0\leq u_d \leq k-1$.
On the other hand, $ ku_k + du_d = \aa_2 <0$ and
$ du_0 + (d-k) u_k = d\aa_1 - \aa_2 <0$, which means that
in either case, $u_k<0$.  Therefore the left-hand side of
(1.16) is never zero and, consequently neither is the
right-hand side.  This implies that $u_0\geq 0$ and
$u_d\geq 0$ and for all $j\geq 0$
$$ c_j = c \, (-1)^{jd} \, \frac {(-u_k+jd)!} 
{(u_0+j(d-k))!\,(u_d+jk)!}
\,,$$
for some constant $c$.    

\noi But a function with such an expansion may not
be rational.
Indeed, Stirling's formula implies  that, asymptotically,
$c_j \sim \mu\, j^{-\aa_1}\,{\lambda^j}/ {\sqrt j}$, 
for some constant
$\mu$. On the other hand, for a rational function 
$R$ whose denominator is a power of
$(1-\lambda\,z)$, we would have $c_j \sim p(j)\, {\lambda^j}$ 
with $p$ a polynomial.
This completes the proof of Theorem~1.9 in the case $m=1$.
\medskip
\noi
Assume now that Theorem~1.9, and its consequence Theorem~1.10, 
are valid
for $m-1$, $m\geq 2$, and consider the $\AA$-hypergeometric 
system associated with
the matrix $\AA$ in (0.4). 

\smallskip 
\noi
Let $R(x)$ be a rational solution with
parameter $\alpha$ in the Euler-Jacobi cone; 
that is, $\aa=(\aa_1,\aa_2)$ with $d\alpha_1<\alpha_2<0$.  
Suppose we can write
$R(x) = {P(x)}/ ({x_d^r\,Q(x))}\,;\  r>0$
and assume that $x_d$ does not divide $P(x)$ or $Q(x)$. 
Let again $x'=(x_0,x_{k_1},\ldots,x_{k_m})$ and
write $Q(x) = \sum_{k=0}^r q_k(x')\,x_d^k$ with
$q_0(x')\not=0$, and set $c_k(x') = - q_k(x')/q_0(x')$. Then
$$R(x) = \frac {P(x)} {x_d^r\cdot q_0(x')} 
\biggl(1+ \sum_{m=1}^\infty\bigl(
 \sum_{k=1}^r c_k(x')\,x_d^k\bigr)^m
\biggr) = \sum_{\ell\geq -r}A_\ell(x')\,x_d^\ell\,.\eqno{(1.17)}$$
By inductive assumption,
since $A_\ell(x')$ is a rational 
$\BB$-hypergeometric function, it must
be a Laurent polynomial as described by Theorem~1.4.
(Note that even though we may have $\ell = \gcd(k_1,\ldots,k_m)>1$ we
may easily reduce to the system associated with a matrix where 
$k_i$ has been replaced by $k_i/\ell$.)
Set $x_d=0$ and 
consider the non-trivial rational function 
$A_{-r}(x') = P(x',0)/Q(x',0)$
which has $\BB$-exponent $\alpha'(r)=( \alpha_1 +r, \alpha_2 + d\,r)$.  
It is clear that $\alpha'(r)\not\in E(\BB)$ 
and hence there may be, up to constant, at most one Laurent
solution with exponent $\alpha'(r)$.  We distinguish the two
possible cases.
 
\medskip
\noi{ i)} Suppose  $A_{-r}(x')$ is a non-zero multiple of
$\Psi^{\BB}_{k_m}(\alpha'(r),x')$.  Then it contains 
  a non-zero term of the form
$$\frac {\hat x^{u} } {x_{k_m}^{s}}\ ; \ \ s>0$$
where $\hat x^u$ is a monomial with positive exponents involving
only the variables $x_0,x_{k_1},\ldots,$ $x_{k_{m-1}}$.
Thus a Laurent series expansion of $R(x)$, as a function of $x_d$, has
a non-zero term of the form
$$\frac {\hat x^{u} } {x_{k_m}^{s}\,x_d^r}\ ;\ \  r>0\,,s>0$$
Successive applications of the fact that
 $R$ is in the kernel of the operator (1.15), with $k_m$ in the
place of $k$, yields that
$R$ must also contain a non-zero term
whose  derivative 
$D_0^{j(d-k_m)}\,D_d^{jk_m}$
is a multiple of
$$\frac {\hat x^{u} } {x_{k_m}^{s+dj}\,x_d^{r}}\ ;\ \  r>0\,,s>0$$
for all $j \geq 0$ which is  impossible as soon as $j\,k_m\geq r$.

\medskip
\noi{ ii)} Suppose  $A_{-r}(x')$ is a non-zero multiple of
$\Psi^{\BB}_{0}(\alpha'(r),x')$.  Then it contains 
  a non-zero term of the form
$$\frac {\check x^{u} } {x_{0}^{s}}\ ;\ \  s>0$$
where $\check x^u$ is a monomial with positive exponents involving
only the variables $x_{k_1},\ldots,x_{k_{m}}$.
which implies that a Laurent series expansion of $R(x)$, 
as a function of $x_d$, has
a non-zero term of the form
$$\frac {\check x^{u} } {x_{0}^{s}\,x_d^r}\ ;\ \  r>0, s>0$$
But, since the operator 
(1.15), with $k_m$ in the
place of $k$ must vanish on $R$,
there must also be a non-zero term of the form:
$$\frac {\check x^{u}\,x_{k_m}^d } {x_{0}^{s+d-k_m}\,x_d^{r+k_m}}\ \ 
$$
which contradicts the index bound in (1.17).
 
\medskip
 
\noi By symmetry we can then assume that $R(x) = P(x)/Q(x)$ and 
neither $x_d$ nor $x_0$ divide $Q$. Thus $R(x)$ is written as
in (1.17) with $r=0$.
For each $\ell\geq 0$, $A_\ell(x')$ is a solution
of the $\BB$-hypergeometric system with exponent $(  \alpha_1 - \ell,
\alpha_2 -d\,\ell)$.  
By inductive hypothesis, these must be Laurent polynomials and
it is easy to check that the only possible denominators are
powers of $x_{k_m}$.  Thus, only powers
of $x_{k_m}$ may appear in the denominator 
of the above expansion for $R$,
which implies that $q_0(x')$ must be of the 
form $x_{k_m}^d$ and, therefore,
$Q(x)$ has bidegree $(d, k_md)$.  On the 
other hand, a symmetric argument
would imply that $Q(x)$ must contain a term 
of the form $x_{k_1}^e$ and
hence $Q(x)$ should have bidegree $(e, k_1 e)$.  
Since $m>1$, this implies
that $Q(x)$ has degree $0$, but then $R(x)$ is a polynomial
solution 
which is impossible since, being in the Euler-Jacobi cone,
$\alpha\not\in \AA\cdot \N^{m+2}$.  \ \ $\diamond$

\bigskip

As we have noted before, given $s\in\Z$, the sum (0.6)
$p_s(x) \ =\ \rho_1^s(x) + \cdots + \rho_d^s(x)$,
 of the powers of the roots of the generic polynomial
(0.5) is a rational $\AA$-hypergeometric function with exponent
$(0,-s)$.  By Theorem~1.10 it must be 
a Laurent polynomial and, therefore expressible in terms of
$\Psi^\AA_0$ and $\Psi^\AA_d$. In fact,
\medskip

\proclaim Corollary 1.12.  
For $s>0$,
$$p_s(x) \ =\ 
s\cdot \Psi_d^\AA((0,-s);x)\ =\ 
s\cdot  \sum_{r=1}^s \ (-1)^r \, (r-1)!\  \frac 
{ \Phi^{\BB}((r,r\,d-s); 
x')} {x_d^r}\,,
\eqno{(1.18)}$$ 
while for $s<0$,
$$p_s(x) \ =\ 
s\cdot \Psi_0^\AA((0,-s);x) \ =\ 
s\cdot  \sum_{r=1}^s \ (-1)^r \, (r-1)!\  \frac { \Phi^{\CC}((r,-s); 
\tilde x)} {x_0^r}\,.\eqno{(1.19)}$$ 

\medskip\noi{\bf Proof:\ } It suffices to consider the normal case,
$d=m+1$, and then set the appropriate variables
equal to zero.
For $s>0$, it follows from (1.11), that $p_s(x)$
must be a multiple of $\Psi_d^\AA((0,-s);x)$.  It is easy to
see that the value of 
the multiple must be
$s$ by specialization to the case when $f(x;t) = t^d + t^{d-1}$.  
The statement for $s<0$ follows by symmetry after observing that
with the change of variables (1.6), the polynomial $\hat f(y;\tau)$
associated with the matrix $\hat \AA$ in (1.5) is given by
$\hat f(y;\tau) = \tau^d\cdot f(x;\tau^{-1})$ and, consequently,
its roots are the inverse of those of $f$.\ \ $\diamond$. 

\medskip
\noi{\bf 1.13. Remarks: } i) Note that each term in the right hand side of
(1.18) is of total degree zero and, therefore, we may express
$p_s(x)$ as a polynomial in
$x_{d-j}/x_d = (-1)^{j}\, \sigma_j(\rho_1,\ldots,\rho_d)$, 
$j=1,\ldots,d$, where
$\sigma_j$ is the $j$-th elementary symmetric polynomial.  This yields
the classical Girard formulas. 

\medskip\noi ii) As we have also noted in the introduction, the total
sum of the local residues  (0.7) gives a rational $\AA$-hypergeometric
function with exponent $(-a,-b)$ and hence, as in Corollary~1.12, 
it must be a multiple of $\Psi_d^\AA((-a,-b);x)$ if $b>0$ and
of $\Psi_0^\AA((-a,-b);x)$ if $b\leq 0$.  In particular, for
$a=1$, $b>0$ we have
$$\sum_\rho {\rm Res}_{\rho(x)}\,\left( \frac {t^b} {f(x;t)}
\frac {dt} {t} \right) \ =\ \cases{
0, &if $0<b<d$;\cr
- \Psi_d^\AA((-1,-b);x), &if $b\geq d$.\cr}$$

\medskip

We end this section with a result that should be seen as a complement
to (1.10) and (1.11) and which will be of use in \S2.
\medskip

\proclaim Proposition 1.14. Let $\aa = (\aa_1,\aa_2) \in \Z^2\setminus 
I(\AA)$ and set $s = s(\aa) := d\aa_1 - \aa_2$.  
\smallskip
\noi i) If $s > 0$, then $\Psi_d^\AA(\alpha;x) = 0$ if and only if
$p_s(x) =0$.
\smallskip
\noi ii) If $\aa_2 > 0$, then $\Psi_0^\AA(\alpha;x) = 0$ 
if and only if
$p_{-\aa_2}(x) =0$.

\medskip\noi{\bf Proof:\ } It suffices to prove the first
statement and then deduce ii) by symmetry.
Suppose $\beta\in \Z^2$ with $\beta_1 = \aa_1 - r$, $r>0$ and
$s(\beta) = s(\alpha)$ then, if
$\Psi_d^\AA(\alpha;x) = 0$  we have by (1.12):
$\Psi_d^\AA(\beta;x) 
= D_d^r \Psi_d^\AA(\alpha;x)
 = 0$.
But, if  $\Psi_d^\AA(\alpha;x)\not = 0$, then the same argument
implies that $\Psi_d^\AA(\beta;x) \not=0$ since the dependence
of $\Psi_d^\AA(\alpha;x)$ on $x_d$ is not polynomial.  Hence, given
$\alpha, \beta \in \Z^2\setminus 
\AA\cdot \N^{m+2}$  such that $s(\alpha) =s(\beta)$, 
 $\Psi_d^\AA(\alpha;x)= 0$ if and only if 
 $\Psi_d^\AA(\beta;x) = 0$.  The result now follows from the
fact that  $\Psi_d^\AA((0,-s);x) = s\,p_s(x)$.\ \ $\diamond$

\bigskip\bigskip

\beginsection {\bf 2. Algebraic  solutions}

\bigskip

In this section we  introduce a family
$\psi_\rho$ of local algebraic solutions
of the $\AA$-hyper\-geometric system associated with a monomial curve.
These solutions, which are given in terms of the roots $\rho(x)$ of the generic
polynomial (0.5), will play a central role in \S 3 when
we compute the
holonomic rank and construct a basis of local solutions for
all  exponents $\aa\in\Z^2$.  
 
\medskip
Let $\AA$ be as in (0.4).
Given an open set $\UU \subseteq \C^{m+2}$, let
 $\HH (\alpha) ({\cal U})$ denote the space of solutions, holomorphic on
$\cal U$, of the $\AA$-hypergeometric system  with exponent $\alpha$, and
$\HH_{alg}(\alpha) (\UU)$  the subspace of algebraic solutions.
We will drop the reference to the open set $\UU$ whenever we are only
interested in the space of local solutions around a generic point.
Let $\rho(x)$ be a root of the polynomial $f(x;t)$ defined
by (0.5), holomorphic for
$x$ in some simply-connected open subset 
$\UU\subset \C^{m+2}\setminus \Sigma$.  

\smallskip

Given 
$\alpha = (\alpha_1,\alpha_2)\in \Z^2$, $\alpha_1\geq 0$ we define:
$$\psi_\rho(\alpha;x) \ :=\  \sum_{{i=0} \atop i\not= 
\alpha_2}^{d\alpha_1}
\Phi^{\AA} ((\alpha_1,i);x) \ { \rho^{i-\alpha_2}(x)  
\over  i-\alpha_2}
 \ +\  \Phi^{\AA} (\alpha;x) \ \log( \rho(x)).
 \eqno{(2.1)}$$
Note that the condition $i \not= \alpha_2$ is automatically
satisfied when $\alpha_2 < 0$ or $\alpha_2 > d\alpha_1$,
and 
that for $\alpha_1 =0$, $\aa_2\not=0$, $\psi_\rho (\alpha)$ is just
$ \rho^{-\alpha_2}  / ( -\alpha_2)$. (If there is no ambiguity
we will drop any reference to the variable $x$.)  

For $\aa\not\in I(\AA)$, the hypergeometric polynomial
$\Phi^{\AA} (\alpha;x)$ vanishes and $\psi_\rho(\alpha)$
is an algebraic function.  This is the case which will be
of interest in this section; however, in \S 3 we will need
(2.1) for arbitrary $\aa\in\Z^2$ and is, therefore, more 
economical to work in this slightly more general setting.
We will extend the definition of 
$\psi_\rho(\alpha)$ to the case $\aa_1 < 0$ in (2.6).

\bigskip
 
\proclaim Proposition 2.1.  Let $\alpha, \alpha' \in \Z^2$ and assume that
$\alpha_1, \alpha'_1 \geq 0$.  Let $u,u'\in \N^{m+2}$ be such that
$\alpha - \AA\cdot u = \alpha'-\AA\cdot u'$.  Then
$$D_u \psi_\rho(\alpha) = D_{u'} \psi_\rho(\alpha')\,.\eqno(2.2)$$
 
\medskip
\noi{\bf Proof:} 
In order to prove (2.2) we show, as a first step, that if 
$\alpha  =(\aa_1,\aa_2)$ and $\alpha_1>0$, then
$$
\frac {\partial} {\partial x_\ell} \bigl(\psi_\rho(\alpha)\bigr) =
\psi_\rho(\aa - \AA\cdot e_\ell). \eqno(2.3)$$
Since, by
(1.2), $\ \partial(\Phi^\AA(\aa))/\partial x_\ell = 
\Phi^\AA(\aa-\AA\cdot e_\ell) $, we have
$$\eqalign{
\frac {\partial} {\partial x_\ell} \bigl(\psi_\rho(\alpha)\bigr) &=
\sum_{{i=0} \atop i\not= \alpha_2}^{d\alpha_1}
\Phi^{\AA} ((\alpha_1-1,i-\ell))
\ { \rho^{i-\alpha_2}  \over  i-\alpha_2}
\  +\ \Phi^\AA(\aa - \AA\cdot e_\ell)\,\log(\rho) 
 \cr
&\quad+\  \sum_{{i=0} \atop i\not= \alpha_2}^{d\alpha_1}
 \Phi^{\AA} ((\alpha_1,i))
  \ { \rho^{i-\alpha_2-1}  }\, \frac {\partial \rho} {\partial x_\ell}
\ + \ \Phi^\AA(\aa)\,\rho^{-1}\,\frac {\partial \rho} {\partial x_\ell}\cr}
\eqno{(2.4)}$$
Note that the last two terms in the expression above cancel since
$f^{\aa_1}(x;\rho(x))=0$ implies that:
$$\sum_{{i=0 \atop \ i \not= \alpha_2}}^{d \alpha_1}
 \Phi^{\AA}(\alpha_1, i) \rho^i\  
=\  - \Phi^{\AA} (\alpha) \rho^{\alpha_2}.$$
On the other hand, setting 
 $j = i-\ell$,
the right-hand side of (2.4) becomes 
 $$
 \sum_{{j=0} \atop j\not= \alpha_2 - \ell}^{d(\alpha_1 - 1)}
 \Phi^{\AA} ((\alpha_1-1,j))
  \ { \rho^{j-(\alpha_2-\ell)}  \over  j-(\alpha_2-\ell)}\  + \ 
\Phi^\AA(\aa - \AA\cdot e_\ell)\,\log(\rho)
\   = \ 
 \psi_\rho(\aa - \AA\cdot e_\ell) 
 \ $$
 where we have used that
 the hypergeometric polynomial $\Phi^{\AA} ((\alpha_1-1,j))$
 vanishes for $j<0$ or $j>d(\alpha_1 - 1)$.  This proves (2.3).

 \medskip\noi
 Applying (2.3) successively we may  assume that
 $\alpha_1 = \alpha_1'=0$. 
Suppose  $\alpha = (0,-s)$, $s\not= 0$,
 so that $\Phi^\AA(\aa)=0$ and
$\psi_\rho(\alpha) = \rho^s/s$.  Given $u\in\N^{m+2}$
 we want to compute $D_u(\psi_\rho(\alpha))$.  Locally on $x$ we
 can write
 $$\eqalign{
 D_u\bigl(\frac {2\pi i} s \,\rho^s\bigr) &=  \int_\Gamma\,
 \frac {t^s} {s} \,D_u\bigl(\frac {f'(x;t)} {f(x;t)}\bigr)\, dt\cr
 &= \int_\Gamma\,
 \frac d {dt} \bigl(\frac {t^s} {s} \,
D_u\bigl(\log {f(x;t)}\bigr)\bigr)\, dt 
 - \int_\Gamma\,
  {t^{s-1}}  \,D_u\bigl(\log {f(x;t)}\bigr)\, dt\cr
 &= (-1)^{\beta_1-1} (\beta_1 - 1)!\,\int_\Gamma\,
 {t^{s-1}}\  \frac {t^{\beta_2}} {f^{\beta_1}}\ dt\ ,
 }\eqno(2.5)
 $$
 where $\AA\cdot u = (\beta_1,\beta_2)$ and
  $\Gamma$ is a sufficiently small loop in the complex
 plane.  Thus $D_u(\psi_\rho(\alpha))$ depends only on
 the pair $(\beta_1, \beta_2 +s) = -(\alpha - \AA\cdot u)$.

\smallskip\noi
If $\alpha=(0,0)$, $\psi_\rho((0,0))=\log(\rho)$ and we can show
that $D_u(\log(\rho))$ depends only on $\AA\cdot u$ arguing as above with
$\log(t)$ taking the place of $t^s/s$.  This 
completes the proof.
$\ \ \diamond$

 \bigskip
 Proposition 2.1 now allows us to extend the definition of 
 $\psi_\rho(\alpha)$ to the case $\alpha_1<0$.  Indeed, let
 $u\in \N^{m+2}$ be such that $\AA\cdot u =(\beta_1,\beta_2)$ with
 $\beta_1\geq -\alpha_1$ and set:
 $$\psi_\rho(\alpha) := D_u(\psi_\rho(\alpha + \AA\cdot u)) \eqno(2.6)$$
 Clearly, this definition does not depend on the choice of $u$ and
for any $\alpha\in\Z^2$, $u\in\N^{m+2}$:
$$D_u(\psi_\rho(\aa))\  =\ \psi_\rho(\aa-\AA\cdot u).\eqno(2.7)$$
If $\aa\not\in I(\AA)$ we can choose $u$ so that
 $\alpha + \AA\cdot u\not\in I(\AA)$ either, and hence
$\psi_\rho(\alpha)$ is an algebraic function.
\bigskip
\proclaim Theorem 2.2.  For $\aa\in \Z^2 \setminus I(\AA)$
the algebraic functions $\psi_\rho(\aa)$ are 
$\AA$-hypergeometric.

\medskip
\noi{\bf Proof:} Given  (2.6), it suffices to consider the
case $\aa_1\geq 0$.
The hypergeometric polynomials 
$\Phi^{\AA} ((\alpha_1,i);x)$ are $\AA$-bihomogeneous of bidegree
$(\alpha_1,i)$ while the powers 
$ \rho^{i-\alpha_2}(x)$ have bidegree $(0,\alpha_2-i)$.  Hence,
$\psi_\rho(\aa)$ satisfies the homogeneity equations (0.2) with
exponent $\aa$.  On the other hand, it is an immediate consequence
of Proposition~2.1 that $\psi_\rho(\aa)$ satisfies the higher-order
equations (0.1).
$\ \ \diamond$

 \bigskip
 
 Let now $\UU\subset \C^{m+2}\setminus \Sigma$ be a simply-connected open
 set and let $\rho_1(x),\ldots,\rho_d(x)$ be distinct
 roots of the polynomial $f(x;t)$, holomorphic for $x\in\UU$.
 Let $\alpha\in \Z^2 \setminus I(\AA)$ and set
 $\psi_j(\alpha) := \psi_{\rho_j}(\alpha) \in \HH_{alg}(\alpha)(\UU)$. 
 The function 
 $$\Psi^\AA(\alpha) := 
 \psi_1(\alpha) + \cdots + \psi_d(\alpha) $$
 is then  a rational solution of the $\AA$-hypergeometric system
 with parameters $\alpha$. It follows from (2.7)
  that, for any $u\in \N^{m+2}$,
 $$
 D_u(\Psi^\AA(\alpha)) = \Psi^\AA(\alpha - \AA\cdot u)\,. \eqno(2.8)$$

 \bigskip
 
 \proclaim Proposition 2.3.  Let 
$\alpha \in \Z^2\setminus I(\AA)$.
 Then,
 \par
 \smallskip\noi
 {i)} $\Psi^\AA(\alpha;x) = \Psi^{\AA}_d(\alpha;x) +
 \Psi^{\AA}_0(\alpha;x),$  defined as in (1.7) and (1.8).
 \par
 \smallskip\noi{ii)} $\Psi^\AA(\alpha) =0$ 
if and only if there exist no
 non-trivial rational $\AA$-hypergeometric functions with parameter
 $\alpha.$
\par
\smallskip\noi{iii)}
If $\alpha =(\alpha_1,\alpha_2)$ and $\alpha_1>0$, then
$$
\Psi^{\AA}_0(\alpha;x)  \ 
=\  \sum_{{i=0} \atop i< \alpha_2}^{d\alpha_1}
\Phi^{\AA} ((\alpha_1,i);x) \,\Psi^{\AA}_0((0,\alpha_2-i);x)\,,
\hbox{\ and}\eqno{(2.9)}$$
$$\Psi^{\AA}_d(\alpha;x)\ 
=\ \sum_{{i=0} \atop i> \alpha_2}^{d\alpha_1}
\Phi^{\AA} ((\alpha_1,i);x) \,\Psi^{\AA}_d((0,\alpha_2-i);x)\,.
\eqno{(2.10)}$$

 \medskip
 
 \noi{\bf Proof:} By Theorem~1.4~iii), $\Psi^\AA(\alpha)$ is a linear
 combination
 $ \lambda\Psi^{\AA}_d(\alpha;x) + \mu\Psi^{\AA}_0(\alpha;x)$. Moreover,
 if $\alpha_1=0$ the result follows from (1.18) and (1.19).  Therefore, 
 computing derivatives with respect to $x_d$,
 the
 result follows for $s(\aa) = d\aa_1 - \aa_2 >0$, 
$\aa_1<0$ where $\Psi_0^\AA(\aa;x)=0$.
 By symmetry i) also holds for $\aa_2>0$, $\aa_1<0$.  
 \smallskip
 \noi If $\aa \in \Z^2\setminus I(\AA)$ is such that 
 $\Psi^{\AA}_d(\alpha;x)\not=0$ then, 
 $$\Psi^\AA((0,-s(\aa))) \ = \ \frac {\partial^{\aa_1} \Psi^\AA(\alpha)} 
 {\partial x_d^{\aa_1}} \ =\  \lambda\, \frac {\partial^{\aa_1} \Psi_d^\AA(\alpha)} 
 {\partial x_d^{\aa_1}}\ =\ \lambda\,\Psi_d^\AA((0,-s(\aa)))$$
 which implies $\lambda=1$ since, because of Proposition~1.14, 
 $\Psi_d^\AA((0,-s(\aa)))\not=0$.  A similar argument shows that if 
 $\Psi^{\AA}_0(\alpha;x)\not=0$ then $\mu =1$.
 \smallskip\noi
 The second assertion is an immediate consequence of i) and iii) in 
 Theorem~1.4. The identities (2.9) and (2.10) follow from part i),
together with (1.18) and (1.19).\ \ $\diamond$

 \bigskip

We now determine the dimension of the subspace of algebraic hypergeometric
functions over $\UU$ spanned by $\psi_1(\aa), \ldots, \psi_d(\aa)$,
$\aa \not \in I(\AA)$.

\smallskip

 \proclaim Theorem 2.4.  For 
$\alpha \in \Z^2\setminus I(\AA)$, the
 $\AA$-hypergeometric functions 
$\psi_1(\alpha),\ldots,\psi_d(\alpha)$
 span a linear space of dimension at least $d-1$.  Moreover, they are
 linearly dependent if and only if $\Psi^\AA(\aa) =\ 0\,.$

\bigskip

Before giving the proof of Theorem 2.4,
we first recall the construction by Gel'fand, Zelevinsky, and
Kapranov [9] 
of $\Gamma$-series solutions for the $\AA$-hypergeometric system
and the expressions obtained by Sturmfels [22] 
for the roots of $f(x;t)$ in terms of them.  We will begin
by considering the normal case, and
eventually we will specialize coefficients to study the
general case.
We only need to consider
the coarsest triangulation of the polytope $P$, the convex
hull of $\AA$ and the origin, i.e. the one consisting of the single
simplex $P$.  As before, we let $\LL$ stand for the integral kernel
of $\AA$, that is, the sublattice of elements $v\in \Z^{d+1}$ such
that $\AA\cdot v = 0$.  Given $u\in\Q^{d+1}$ we define the
formal power series
$$
\bigl[ x_0^{u_0} x_{1}^{u_{1}} \cdots x_d^{u_d} \bigr] \quad := \quad
\sum_{v \in { \LL}}
\prod_{i=0}^d  \, \bigl(\,\gamma(u_i,v_i) \,x_i^{u_i+v_i} \bigr)\,, 
\eqno (2.11) $$
where,
for any rational number $u$ and any integer $v$, we write
$$ \gamma(u,v) \,\,\, := \,\,\, \cases{ \qquad
   1         &  if $v = 0$, \cr
   u (u  -  1) (u  - 2)  \cdots  (u+v+1)    & if $v  < 0$, \cr
   0          & if $u$ is a negative integer and 
$u \geq -v $, \cr \qquad
  {1 \over (u+1)(u+2) \cdots (u+v)} & otherwise. \cr }$$
If $u$ has no negative integer coordinates, or is of the form
$(0,\ldots,0,1,-1,0,\ldots,0)$,  
the series $\bigl[ x_0^{u_0}\cdots x_d^{u_d} \bigr]$ 
is a formal solution of the  ${\AA}$-hypergeometric system
with parameters ${ \AA} \cdot u \in \Q^2$ 
(see [9,~Lemma 1], [22,~Lemma 3.1]).  
Moreover, if at most two of the exponents $u_i$ are non-integers,
the series (2.11) converges in a suitable open subset of
$\C^{d+1}\setminus\Sigma$ ([9]). 
\medskip
Let $\xi_1,\ldots,\xi_d$ be the $d$-th roots of $-1$
and consider the series
$\  \rho_{i}(x)\,:= \,\sum_{a=1}^d\,\xi_i^a\cdot \sigma_a(x)\,,$
where
$$
\sigma_1(x) \ := 
\ \biggl[ {x_{{0}}^{1 \over d}  x_{d}^{-{1 \over d}}} 
 \biggr]  ; \quad\quad
 \sigma_a(x) \ := \ {1 \over d} \cdot \biggl[
 x_{a-1} \,x_{0}^{{a-d \over d}}  x_{d}^{-{a \over d}}
 \biggl],\ a=2,\ldots,d\,.$$
It follows from [9,~Proposition 2] that there exists an
open set  $\VV\subset \C^{d+1}\setminus \Sigma$ of the form
$$\VV \ :=\ \{x\in \C^{d+1} \,:\, |x_0|^{d-j} \, 
|x_d|^{j} > M\,|x_j|^{d}\ ;
 \ j \not= 0,d\}, $$
for some positive real constant $M$, where all these
series converge, and, according to [22,~Theorem 3.2], locally on $\VV$, 
they define the holomorphic 
$d$ roots of the generic polynomial
$\sum_{j=0}^d \,x_j\cdot t^j$.  Given a positive integer $s$ we
consider the powers $\rho_i^s(x)$ and write
$$
  \rho_{i}^s(x)\,= \,\sum_{b=1}^d\,\xi_i^b\cdot \theta_b(x)\quad
 \hbox{where}\quad
 \theta_b(x) = \sum_{a_1+\cdots+a_s = b + \ell d} \,(-1)^\ell\,
 \prod_{j=1}^s \sigma_{a_j}(x)$$
 
\bigskip
We now consider the $\AA$-hypergeometric system associated
with the matrix (0.4) and recall  that we are
assuming that $\gcd(k_1,\ldots,k_m,d)=1$.
Let $J$ denote the complement of $\{0,k_1,\ldots,k_m,d\}$ in 
$\{0,1,\ldots,d\}$, and $V_J$ the $(m+2)$-dimensional
subspace of $\C^{d+1}$ obtained by setting $x_j=0$, $j\in J$.
Note that $\VV \cap V_J $ is non-empty.
 
\bigskip
\proclaim Lemma 2.5.  
For $b=1,\ldots,d-1$, the restriction  of 
$\theta_b(x)$ to $V_J$ is non-trivial.
The same is true of any of its derivatives 
${D_u \theta_b }$,
$u\in\N^{m+2}$, with respect to variables $x_i$, $i\not\in J$.

\bigskip\noi
{\bf Proof:\ } Note that if 
$\rho^s = (\rho_1^s,\ldots,\rho_d^s)^T$ and 
$\theta = (\theta_1,\ldots,\theta_d)^T$, then
$\theta = M^{-1}\cdot  \rho^s$ where $M$ is the non-singular matrix
$M=(\xi_i^a)$, $i,a=1,\ldots,d$. In particular, 
$\theta_b(x) \in  \HH ((0,-s)) (\UU)$.  
 
\medskip\noi
We now claim that for any set of indices 
$a_1,\ldots,a_s$ 
such that $a_1+\cdots+a_s = b + \ell d$, for some
$\ell\in \N$,
$$\theta_b(x) = \lambda\cdot \biggl[
x_{a_1-1} \,\cdots x_{a_s-1} \,x_{0}^{{b+(\ell-s)d \over d}}  
 x_{d}^{-{b+\ell d\over d}}
 \biggl]\,,\eqno(2.12)$$
for some non-zero constant $\lambda$.  Indeed, forgetting for
the moment the coefficients, suppose that
$x^w = x^{w^{(1)}}\cdots x^{w^{(s)}} $ 
is a monomial appearing in the product
$\prod_{j=1}^s \sigma_{a'_j}(x)$ where $a'_1 + \cdots + a'_s = 
b + \ell' d$.  Then $\AA\cdot w = (0,-s)^T$ and therefore, 
$x^w$ must differ from the monomial inside the bracket (2.12) by
a monomial of the form $x^v$ with $v\in \LL$.  
This means that all the
monomials in the power series of $\theta_b(x)$ 
appear in the $\Gamma$-series
of (2.12).  But, on the other hand, since we already know that
$\theta_b(x)$ is $\AA$-hypergeometric, if a monomial, such as the
one in (2.12), 
appears in its expansion then the whole $\Gamma$-series
must appear and with the appropriate coefficients. 
 
\medskip\noi
Suppose now that we set the variables $x_j$, $j\in J$, equal to zero.
We may assume without loss of generality that in the 
bracket in (2.12),
$a_1-1 ,\ldots, {a_r-1}$ are the only indices in $J$.  Since 
$\gcd(k_1,\ldots,k_m,d)=1$, there exist positive integers 
$p'_d,p_d'',p_1,\ldots,p_m$ such that
$$(a_1 -1)+\cdots + ({a_r} - 1) + p'_d\,d\ 
 =\  p_1\,k_1 + \cdots + p_m\,k_m + p_d''\,d$$
Setting $p_d = p_d'' - p_d'$, there exists $p_0\in \Z$ such that
$$v \ :=\ p_0 \,e_0 +  p_1\,e_{k_1} + \cdots 
+ p_m\,e_{k_m} + p_d\,e_d - 
 e_{a_1 -1}-\cdots - e_{a_r -1} \ \in\  \LL\,.$$
 Consequently, multiplication of the monomial in the bracket
(2.12) by $x^v$ yields 
 a term in the $\Gamma$-series 
 which does not involve any variables from the
 index set $J$.  On the other hand, it is easy to check that
 all coefficients $\gamma(u_i,v_i)$ are non-zero and therefore the
 restriction of $\theta_b(x)$ is non-trivial.
 
 \medskip\noi
 The statement about the derivatives of $\theta_b(x)$ follows
  from the fact that for $b<d$, the exponents of
 $x_0$ and $x_n$ in the bracket in (2.12) are not integers, while the
 exponents with which any of the other variables $x_{k_i}$ appears
 in the $\Gamma$-series cannot be bounded since for any $\ell\in \N$,
 the element
 $$ \ell\cdot\bigl(d\,e_{k_i} - (d-k_i)\, 
e_0 - k_i\,e_d\bigr) \ \in \ \LL\,.
 \eqno{\diamond}$$

 \bigskip\noi{\bf Proof of Theorem 2.4:  \  }We consider, first of all, the case
 $\alpha = (0,-s)$ with $s$ a positive integer.  Then 
 $\psi_j(\alpha) = \rho_j^s/s$ and, since 
 $\rho^s = M\cdot \theta$ with $M$ non singular,
$\rho_1^s,\ldots,\rho_d^s$
 are linearly independent if and only if 
$\theta_1(x),\ldots,\theta_d(x)$
 are so.  But comparing the exponents of
 $x_d$ in the corresponding
 $\Gamma$-series it is clear that
 $\theta_1(x),\ldots,\theta_d(x)$ will
 be linearly independent unless their
 restriction to $V_J$ vanishes.  On the other hand, 
 Lemma~2.5 asserts that only $\theta_d(x)$ may be
identically zero when
restricted to $V_J$.  Hence, the dimension of the linear span
of $\rho_1^s(x),\ldots,\rho_d^s(x)$ is at least $d-1$ and will be
exactly $d-1$ if and only if $\theta_d(x)=0$. But,
$$\sum_{i=1}^d \,\rho_i^s(x) \ 
=\ \sum_{b=1}^d\bigl( \sum_{i=1}^d \xi_i^b
\bigr)\cdot \theta_b(x) \  =\ -d\,\theta_d(x)\,.$$
Thus, $\theta_d(x)=0$ if and only if  
$\rho_1^s(x)+\cdots+\rho_d^s(x)\ =\ 0 $.  
 
\medskip\noi  The assertion for $\alpha = (0,s)$
with $s$ a positive integer
follows from symmetry.  In view of the
definition (2.6), the statement for $\alpha_1<0$
follows from that for $\alpha_1=0$ using the assertion
in Lemma~2.5 about
the derivatives of the $\Gamma$-series $\theta_b(x)$.
 
\medskip\noi  It remains to consider the case $\alpha_1>0$.
Suppose there is
a non-trivial linear relation
$\sum_{i=1}^d \lambda_i\cdot\psi_i(\alpha) = 0 $; because
of (2.7), applying the 
derivative $\partial^{\aa_1}/\partial x_0^{\aa_1}$ we obtain
$$\sum_{i=1}^d \,\lambda_i\cdot\psi_i((0,\alpha_2)) \ =\  0 \,.$$
But $\alpha\not\in I(\AA)$ implies that 
$(0,\alpha_2)\not\in I(\AA)$, 
and therefore, $\lambda_1=\cdots=\lambda_d$,
and the proof is complete.\ \ $\diamond$

\bigskip

The following result implies that the holonomic rank of the
$\AA$-hypergeometric system is at least $d+1$ for all
exponents $\aa \in E(\AA).$

\medskip
\proclaim Theorem 2.6. For any $\aa \in E(\AA), \ 
\dim(\HH_{alg}(\aa)(\UU)) \geq d+1.$

\medskip\noi{\bf Proof: \ }
For any $\aa \in E(\AA)$ the $\AA$-hypergeometric Laurent 
polynomials $\Psi_0^\AA(x)$ and $\Psi_d^\AA(x)$ are
 both non-trivial.  In particular, $0<\aa_2<d\,\aa_1$. From i) in 
Proposition~2.3 we have that
 $$\psi_1(\alpha) + \cdots + \psi_d(\alpha)\  \not=\ 0\,,$$
 and Theorem~2.4 implies that 
 $\psi_1(\alpha),\ldots , \psi_d(\alpha)$ are linearly
 independent.  Moreover, we will show next that so are the functions
 $\psi_1(\alpha),\ldots , \psi_d(\alpha),\Psi^\AA_0(\alpha)$. 
Suppose there is a non-trivial 
linear combination
 $\ \Psi^\AA_0(\alpha) = \sum_{i=1}^d \lambda_i\cdot \psi_i(\alpha)$,
 $\,\lambda_i\in\C$.  By differentiation we obtain a similar relation
$$\Psi^\AA_0((0,\alpha_2))\  = \ 
 \frac {\partial^{\aa_1} \Psi^\AA_0(\alpha)} {\partial x_0^{\aa_1}}
 \  = \  \sum_{i=1}^d \lambda_i\cdot 
 \frac {\partial^{\aa_1} \psi_i(\alpha)} {\partial x_0^{\aa_1}}
 \  = \  \sum_{i=1}^d \lambda_i\cdot \psi_i((0,\aa_2))\,.$$
But (1.10) implies that $\Psi_d^\AA((0,\aa_2))=0$
and therefore $\Psi_0^\AA((0,\aa_2)) =  
\sum_{i=1}^d\psi_i((0,\aa_2))$.  
Since according to
 ii) in Proposition~1.14,
 $\Psi^\AA_0((0,\alpha_2))\not=0$, we have that all
 $\lambda_i=1$.  This implies that 
$\Psi_0^\AA(\aa) = \sum_{i=1}^d\psi_i(\aa) = \Psi_0^\AA(\aa) +
\Psi_d^\AA(\aa)$, i.e. $\Psi_d^\AA(\aa)=0$,
which contradicts our assumption. \ \ $\diamond$
 
\bigskip

\proclaim Corollary 2.7.  The toric ring 
$R = \C[\xi_0,\ldots,\xi_d]/\II_\AA$
is Cohen-Macaulay if and only if, for every $\alpha\in\C^2$
the dimension of the space of $\AA$-hypergeometric functions of
exponent $\aa$, at a generic point, is equal to $d$.  
 
\medskip\noi{\bf Proof:\ } The only if direction is Theorem 2 in 
[9] (see also [10] and [2]).  To prove the converse we note that,
because of Proposition~1.6, if $R$ is not Cohen-Macaulay then
$E(\AA) \not= \emptyset$, and the result follows from Theorem
2.6.

\bigskip\bigskip
\beginsection{\bf 3. Bases of solutions and holonomic rank}

\bigskip
In this section we will  determine the holonomic
rank of the $\AA$-hypergeometric system associated with
a monomial curve for all integral exponents and exhibit
explicit bases of local solutions constructed in terms of
the roots of the generic polynomial (0.5).  
\bigskip
Four different scenarios need to be considered:
\smallskip
\item {$\bullet$} The exponent $\aa\in I(\AA)$:  In this 
case $r(\aa) = d$ and we construct in Theorem~3.1,
$d-1$ local solutions which, together with the 
hypergeometric polynomial $\Phi^\AA(\aa)$ define
a basis of solutions.
\smallskip
\item {$\bullet$} $\aa\in E(\AA)$:  We now have 
$r(\aa) = d+1$ and we have constructed $d+1$
(algebraic) local solutions in Theorem~2.6.
\smallskip
\item {$\bullet$} $\aa\in (E_0(\AA) \cup E_d(\AA))\setminus E(\AA)$:
The holonomic rank equals $d$ and we have from 
Theorem~2.4 a basis of algebraic solutions.
\smallskip
\item {$\bullet$} $\aa\in J(\AA)$: Then $\RR(\aa)=\{0\}$, 
$r(\aa)=d$, and we
construct a basis of local solutions in Theorem~3.5.
\bigskip
In Theorem~3.7 we determine
the holonomic rank $r(\aa)$ for all $\aa\in\Z^2$.
Our starting point is a result of
Adolphson ([2,Corollary~5.20]) which states that even
without assuming 
that the ring $R$ is Cohen-Macaulay, the equality
$r(\aa) = {\rm vol}(P)$ holds for so-called
{\sl semi-nonresonant} exponents $\aa$.  In our particular
case, this condition is equivalent to 
$\aa$ being  in the Euler-Jacobi cone (1.13).

\bigskip

We consider first the case when $\alpha\in I(\AA)$, i.e. when
the $\AA$-hypergeometric polynomial $\Phi^\AA(\aa;x) \not= 0$. 
By Proposition~2.1, given a root $\rho(x)$, the function
$\psi_\rho(\aa;x)$
satisfies the higher-order equations (0.1) but, clearly,
 not the homogeneity
equations (0.2).  However, if we fix a choice of a root 
$\hat\rho$ then,
for any other root $\rho$, the function
$$\tau_\rho(\aa)\ :=\ 
\sum_{{i=0} \atop i\not= \alpha_2}^{d\alpha_1}
 \Phi^{\AA} (\alpha_1,i) \ \frac {\rho^{i-\alpha_2}
- \hat\rho^{i-\alpha_2} } { i-\alpha_2}
 \ +\  \Phi^{\AA} (\alpha) \ \log( \rho/\hat\rho).
 \eqno{(3.1)}$$
is $\AA$-hypergeometric with exponent $\aa$. Indeed, it differs from
$\psi_\rho(\aa) - \psi_{\hat\rho}(\aa)$ by a constant multiple
of $\Phi^\AA(\aa)$.

\proclaim Theorem 3.1. Given $\alpha \in I(\AA)$ and  
a choice of a root $\ \hat{\rho}\ $ of $f(x;t)$ on $\UU$, 
the functions  $\Phi^{\AA}(\alpha) $ and $\tau_\rho(\alpha) $ 
where $\rho$ runs over all roots of $f(x;t)$ on $\UU$ different 
from $\hat{\rho}$, are linearly independent $\AA$-hypergeometric 
functions. 
 
\medskip
 
\noi{\bf Proof:} Suppose $\aa=\AA\cdot u, u \in \N^{m+2},$ 
and suppose there is a 
non-trivial linear combination 
$$\lambda\,\Phi^\AA(\aa) \ +\  
\sum_{\rho\not=\hat\rho} \lambda_\rho \, \tau_\rho(\aa)\ =\ 0.$$
Applying the operator $D_u$ and using (1.2) and (2.7), we may 
assume that $\aa=(0,0)$ and, consequently, $\Phi^\AA(\aa) =1$,
$\tau_\rho(\aa) = \log(\rho/\hat\rho)$, and
$$
\lambda \ +\  
\sum_{\rho\not=\hat\rho} \lambda_\rho \, \log(\rho/\hat\rho)\ =\ 0.
\eqno{(3.2)}$$
Implicit differentiation of the equation $f(x;\rho(x)) = 0 $ yields
${\partial \rho}/{\partial x_\ell} = -\rho^\ell/f'(\rho)$, 
$\ell\geq 0$ and,
consequently,
$$\frac {\partial \log \rho} {\partial x_d}\  
= \ \frac 1 d \frac {\partial \rho^d} {\partial x_0}\ =\ 
\psi_\rho((-1,-d)).$$  
Hence, derivating (3.2), we obtain
$$
\sum_{\rho\not=\hat\rho} \lambda_\rho \,(\psi_\rho((-1,-d))-\psi_
{\hat\rho}((-1,-d)))\ =\ 0,
$$
which, in view of Theorem~2.4, implies $\lambda_\rho=0$
 for all $\rho\not=
\hat\rho$.\ \ $\diamond$

\bigskip

Suppose now that $\aa\in J(\AA)$.
In particular, $\aa\not\in I(\AA)$,
and, by ii) in Proposition~2.3, $\Psi^\AA(\aa)=0$. 
Recall also that  this case includes all  integral
exponents lying in the Euler-Jacobi cone.

As before, we let $\UU$ denote a simply-connected open set
in $\C^{m+2} \setminus \Sigma$ and let $\rho_1,\ldots,\rho_d$
denote the roots of $f(x;t)$ for 
$x\in\UU$.
 Given  $\aa\in \Z^2$ such that
$\aa_1\geq 0$,
we define
$$\chi(\aa)\ :=\ \sum_{j=1}^d \psi_j(\aa)\,\log(\rho_j)\,.
\eqno{(3.3)}$$

\bigskip
\proclaim Proposition 3.2.  Suppose $\aa\in J(\AA)$ is such that
$\aa_1\geq 0$.  Then
the function
$\chi(\aa)$ is $\AA$-hypergeometric with exponent $\aa$.

 \medskip
 \noi{\bf Proof:}  Since 
 $\Psi^\AA(\aa) = \sum_{j=1}^d\,\psi_j(\aa) =0$ it follows that
$\chi(\aa)$ satisfies the equations (0.2) with exponent $\aa$.
In order to check that the higher-order equations (0.1) are
satisfied as well we show, first of all, that if $\aa_1>0$, then
$$\frac {\partial \chi(\aa)} {\partial x_\ell}  = 
\chi(\aa - \AA\cdot e_\ell).\eqno{(3.4)}$$
Indeed,
 $$\eqalign{
 \frac {\partial \chi(\aa)} {\partial x_\ell} \  &= 
\  \sum_{j=1}^d \psi_j (\alpha - \AA\cdot e_\ell) \log(\rho_j)
\  +\  \sum_{j=1}^d \psi_j(\alpha)\, \rho_j^{-1}\,
\frac {\partial \rho_j} {\partial x_\ell}  \cr
 {} & = \ \chi(\alpha - \AA\cdot e_\ell)\  +\  \sum_{i=0}^{d \alpha_1}
 {\Phi^{\AA}((\alpha_1,i)) \over i-\alpha_2}\, 
 \left(\sum_{j=1}^d \rho_j^{i-\alpha_2 - 1} 
\frac {\partial \rho_j} {\partial x_\ell} \right).}$$
 We claim that the second summand is identically zero.
 In fact,
 $$\sum_{j=1}^d \rho_j^{i-\alpha_2 - 1} 
\frac {\partial \rho_j} {\partial x_\ell}  =
\frac {\partial } {\partial x_\ell}  
\left( \sum_{j=1}^d{ \rho_j^{i-\alpha_2} \over i-\alpha_2} \right).$$
Assume $\Phi^{\AA}(\alpha_1,i) \not= 0,$ i.e.,
there exists $w\in \N^{m+2}$
such that  $\AA\cdot w = (\alpha_1,i).$ If $ \sum_{j=1}^d \rho_j^{i-\alpha_2}$
does not vanish identically,  
there exists $v \in \Z^{m+2}$ such that $v_{k_i} \geq 0$
 for all $i=1,\ldots, m$, and either
$v_0\geq 0$ or $v_d\geq 0$, verifying $\AA\cdot v = (0, \alpha_2-i).$
 But then $\AA\cdot (v+w) = \alpha,$ which contradicts the fact that
 $\Psi^{\AA}(\alpha) = 0$. This proves (3.4).
\smallskip\noi
Arguing as in the proof of  
Proposition~2.1, the proof of Proposition~3.2 now reduces
to the following:
\bigskip
\proclaim Lemma~3.3. Given $s\in \Z$, $s\not=0$, and $u\in\N^{m+2}$,
 let
$\beta = \AA\cdot u$ and
$\gamma = (0,-s) - \beta$.   Then, if $\Psi^{\AA}(\gamma) =0$,
the derivative
$$D_u\bigl(\sum_{j=1}^d\, \frac {\rho_j^s} {s} \,\log\rho_j \bigr)$$ 
depends only on $\gamma$.

\medskip\noi
{\bf Proof: } We argue as in (2.5); locally on $x$,
$$\eqalignno{
2\pi i\,D_u\bigl(\sum_{j=1}^d\, 
\frac {\rho_j^s} {s} \,\log\rho_j \bigr)\ &=
\sum_{j=1}^d\int_{\Gamma_j} 
 \frac {t^s} {s} \,\log t \ D_u\bigl(\frac {f'(x;t)} {f(x;t)}\bigr)\, 
dt\cr
&=- \sum_{j=1}^d\int_{\Gamma_j}
\frac d {dt} \bigl(\frac {t^s} {s}\log t\bigr)\ D_u(\log f(x;t))\,dt
&(3.5)\cr
&= c \cdot \left[ \sum_{j=1}^d
\int_{\Gamma_j} 
\log t \ \frac {t^{\beta_2 + s -1}} {f^{\beta_1}} dt
+ \sum_{j=1}^d 
\int_{\Gamma_j} \frac {t^{s-1}} {s} 
\ \frac {t^{\beta_2}} {f^{\beta_1}}dt \right],\cr}
$$
where $\Gamma_j$ is a sufficiently small loop in the complex plane
enclosing only the root $\rho_j$ and $c= (-1)^{\beta_1 -1} (\beta_1 -1)!$.  
Now, according to (2.5), the 
last sum in (3.5), agrees, up to constant, with
$$\sum_{j=1}^d D_u(\rho_j^s/s)\  =
 \  \sum_{j=1}^d D_u(\psi_j((0,-s)))\  = \ 
\sum_{j=1}^d \psi_j(\gamma)\  = \ \Psi^\AA(\gamma)\ =\ 0\,.$$
Hence the Lemma, and Proposition~3.2, follow.\ \ $\diamond$

\bigskip

We now extend the definition (3.3) of $\chi(\aa)$
to the case $\aa\in\Z^2$, $\aa_1<0$, by setting:
$${\chi(\aa) \ =\ \cases{
D_0^{-\aa_1}\,\chi((0,\aa_2)),&if $\aa_2\not=0$\cr
&\cr
D_{k_1}D_0^{-\aa_1-1}\,\chi((0,k_1)),&if $\aa_2=0$,\cr}}\eqno(3.6)$$
where, we recall
$$\chi((0,s)) = \sum_{j=1}^d \frac{\rho_j^{-s}}{-s} \log \rho_j\ ;\quad
s\not = 0\,. \eqno{(3.7)}$$

\bigskip
\proclaim Proposition 3.4. Let $\aa\in \Z^2$, $\aa_1<0$. Then
 \par
 \smallskip\noi
 {i)}  $\ \chi(\aa)(t*x)\ =\ 
t^{\aa_2}\, \chi(\aa)(x) - t^{\aa_2} \,\log t \,\Psi^\AA(\aa;x) $, where
$t\in\C^*$, and $$t*x=
(x_0,t^{k_1}x_{k_1},\ldots,t^{d}x_d)\,.$$
 \par
 \smallskip\noi
 {ii)} \ If $\aa\in J(\AA)$, $\chi(\aa)$ is $\AA$-hyper\-geometric
with exponent $\aa$.
\par
 \smallskip\noi
{iii)}\ For $M$ sufficiently large and $j=1,\ldots,m$,
$$D_0^{M(d-k_j)} D_d^{Mk_j} (\chi(\aa))= D_{k_j}^{Md}
(\chi(\aa))\,.
$$

\medskip\noi{\bf Proof: } In view of (3.7), 
the first assertion follows from (3.6)
together with the fact   that $\rho_j(t*x)
= t^{-1}\,\rho_j(x)$, $j=1,\ldots,d$.
\smallskip\noi If $\aa\in J(\AA)$, $\Psi^\AA(\aa;x)=0$ and it follows
from i) that $\chi(\aa)$ satisfies the homogeneity equations (0.2).
On the other hand, if $\aa\in J(\AA)$, the same is true of 
$\aa - \AA\cdot u$ for every $u\in \N^{m+2}$.  Hence it follows
from Lemma~3.3 that $\chi(\aa)$ satisfies the equations (0.1).
\smallskip\noi The last assertion follows again from Lemma~3.3
for $M$ so that 
$\aa - (M\,d,M\,d\,k_j) \in \EE\JJ \subset J(\AA)$.
\ \ $\diamond$

\bigskip
\proclaim Theorem 3.5. Let $\aa\in J(\AA)$. 
Let $\UU$, $\rho_1,\ldots,\rho_d$, be
as above. Then the functions $\psi_1(\aa),\ldots,\psi_{d-1}(\aa),\chi(\aa) \in
\HH(\aa)(\UU)$ are
linearly independent. 

\medskip\noi {\bf Proof: }
In view of Theorem 2.4, it will be enough to show that $\chi(\aa)$ is
not an algebraic function. In fact, we will show that its orbit
under the monodromy action of $\pi_1 (\C^{m+2} \setminus \Sigma)$ is
infinite. 

\smallskip\noi
For generic values $a_{k_1},\ldots, {a}_{d}$ the polynomial
$f((0, a_{k_1},\ldots, {a}_{d}); t)$ will have a 
a root of multiplicity
$k_1$ at the origin and  $d-k_1$ simple, non-zero roots.  Thus,
for $|a_0|$ small, relative to $|a_{k_1}|,\ldots, |{a}_{d}|$,
the polynomial $f(a;t)$ will have simple roots and $k_1$ of them,
say $\rho_1,\ldots,\rho_{k_1}$,
will be very close to the
$k_1$-th roots of $- {a_0}/
{a}_{k_1}$.  This means that analytic continuation
of the roots along the loop 
$$\gamma(\theta) := ( \exp({2\pi i k_1 \theta}) {a_0}, {a}_{k_1},
\ldots,{a}_d) \ ; \quad \theta \in [0,1]$$ 
returns the roots to their original values, but, for any choice of 
logarithm for all roots of $f$ near $a$:

$$\gamma^*(\log(\rho_j)) = \cases{\log(\rho_j) + 2 \pi i,
&if $j=1,\ldots,k_1$;\cr
\log(\rho_j), &otherwise.\cr}$$
Since $\aa\not\in I(\AA)$, it follows from  (2.1) and (2.6)
that $\gamma^*(\psi_j(\aa)) = \psi_j(\aa)$ for any  
$j =1,\dots,d .$ Therefore, for $\aa_1 \geq 0$,
$$\gamma^*( \chi(\aa)) =
\gamma^*\left(\sum_{j=1}^d \psi_j(\aa)\,\log(\rho_j)\right) =
\chi(\aa) + 2 \pi i \sum_{j=1}^{k_1} \psi_j.$$
Since $0 < k_1 < d$, by Theorem 2.4  $\psi_1 + \ldots
+ \psi_{k_1} \not= 0 \ $, and therefore the orbit of $\chi(\aa)$
under successive powers of $\gamma$ is infinite.

\smallskip\noi
Suppose now that $\aa_1 <0$. Then,
$$ \chi(\aa) =
\frac {\partial^{-\aa_1}} {\partial x_0^{-\aa_1}}
\left(\sum_{j=1}^d \psi_j((0,\aa_2)) \,\log \rho_j \right)
= \sum_{j=1}^d \psi_j(\aa) \, \log \rho_j + R(\aa) .$$
It is straightforward to check that $R(\aa)$ is algebraic
and invariant under the monodromy action, i.e. $R(\aa)$ is
a rational function. Since we have just shown that the
function
$ \sum_{j=1}^d \psi_j(\aa) \, \log \rho_j  \ $
is not algebraic, the proof is complete.
\ \ $\diamond$ 

\bigskip

\noi{\bf 3.6. Remark:} 
Note that since the function $\chi(\aa)$ is not algebraic, 
any rational
$\AA$-hyper-geometric  function $R$ with exponent
$\aa$ in the Euler--Jacobi cone must be a linear combination
of $\psi_1(\aa),\ldots,\psi_d(\aa)$.  On the other hand,
with similar arguments as in the proof of
Theorem~3.5, it is possible to show the existence of a
loop $\gamma$ whose action on the roots is a cyclic permutation.
It is then easy to see that $R$ must be a constant multiple of 
$\sum_{j=1}^d \psi_j (\aa)$, and therefore
it must vanish.  This gives an alternative proof of Theorem~1.9.

\bigskip

\proclaim Theorem 3.7.   For every $\aa\in\Z^2$,
$$d \leq r(\aa) \leq d+1\,.$$
Moreover, $r(\aa)=d+1$ if and only if $\aa \in E(\AA)$.

\medskip
\noi{\bf Proof:} Note, first of all, that the lower
bound follows from
Theorem~2.4 (for $\aa\in E_0(\AA) \cup E_d(\AA)$),
 Theorem 3.1  (for $\aa \in I(\AA)$),
 and Theorem~3.5 (for $\aa \in J(\AA)$).

\smallskip\noi
Suppose now that
$\aa$ is in the Euler-Jacobi cone
$\EE \JJ$. Then, as we observed before, $\aa$ is semi-nonresonant 
in the sense of Adolphson and it follows from 
[2,~Corollary~5.20] that $r(\aa) =d$.
For any $\aa \in \Z^2$, there exist $u \in \N^{m+2}$ such
that $\aa - \AA . u$ lies in $\EE \JJ$, and, because
of Theorem~1.9, for any such $u$
the kernel of the linear map
$$D_u : \HH(\aa) \rightarrow \HH(\aa - \AA.u),$$
contains $\RR(\aa)$.
We will  determine the dimension of $\HH(\aa)$
by computing the kernel and the image of $D_u$ for suitable
$u$.

\smallskip\noi
Suppose first that $\aa\in J(\AA)$. 
For $u = \ell\,e_{k_1}$, $\ell>>0$,
we have $\aa - \AA.u \in \EE \JJ$ and, it follows
from Proposition~1.2 and Corollary~1.5 that
${\rm ker}(D_u) = \RR(\aa) = \{0\}$.  
Therefore, 
$D_u$ is a monomorphism, which implies that $\dim \HH(\aa) 
\leq d$. Since it is at least $d$, we deduce $\dim \HH(\aa) 
=d$. 

\smallskip\noi
Suppose now that $\aa\not\in J(\AA)$, then
$\dim(\RR(\aa)) = 1$  or  $2$.
We can again choose $u = \ell\,e_{k_1}$, $\ell>>0$, so
that $\beta:=\aa - \AA.u \in \EE \JJ$.
As the kernel of $D_u$ is precisely $\RR(\aa)$
it will be enough to show that, for some $\ell$
sufficiently large, the
dimension of the image of $D_u$ is $d-1$.

\smallskip\noi
From (2.6) and 
Theorem~3.5, we deduce that
the functions $\psi_j(\bb)$ generate a subspace of the
image of dimension
at least
$d-1$. The proof will be completed
by showing that the function  $\chi(\beta)$
defined in  (3.6) is not in the image
$D_u(\HH(\aa) )$.

\smallskip\noi
Consider first the case $\aa\in E_0(\AA)$.
Choosing $\ell = s\,d$, we factor
$D_u = D^{sd}_{k_1} = D_d^{ k_1s}  D_0^{ (d-k_1)s}.$
Set $\aa':= \aa - ((d-k_1)s,0)$. It is  enough to
show that $\chi(\beta) \not\in D_d^{ k_1s}(\HH(\aa') )$.  Note
that $\aa'\in E_0(\AA)$ as well and, therefore,
we may assume without loss of generality that $\aa_1 < 0$
and $\beta = \aa - \ell e_d \in \EE \JJ$ for
some sufficiently big $\ell$.

\smallskip\noi
Let $\chi(\aa)$ be as in (3.6).  Since $\beta\in \EE\JJ$, it 
follows from Lemma~3.3 that
$D_d^{\ell}(\chi(\aa)) =
\chi(\beta)$.
Therefore, if $\phi\in\HH(\aa) $ is such that
$D_d^\ell(\phi) = \chi(\beta)$ we must have
$$\phi =  \chi(\aa) + F$$
where $F$ depends polynomially on $x_d$.
On the other hand, because of iii) in
Proposition~3.4 and the fact that $\phi$ is
hypergeometric, we have
$$D_{k_j}^{Md}(F) = D_{0}^{M(d-k_j)}D_{d}^{Mk_j}(F)$$
for all $M$ large enough.  This implies that $F$ is polynomial
on $x_{k_1},\ldots,x_{k_m}$ as well.

\smallskip
\noi But, it follows from i) in Proposition~3.4 that
$$F(t*x) = t^{\aa_2}\,F(x) +
 t^{\aa_2} \log t \,\Psi(\aa)(x)$$
which is impossible since the fact that the action of $t$
does not affect $x_0$ implies that
$F(t*x)$ is polynomial in $t$.

\smallskip
\noi
By symmetry, the result also holds for $\aa\in E_d(\AA)$.
Thus, it remains to consider the case
$\aa  \in I(\AA)$.  For $\ell$ large enough, so that
$\aa'_1 = \aa_1 - \ell(d-k_1) < 0$, we have
$\aa' = (\aa'_1,\aa_2) \in E_0(\AA)$
and an argument similar to the one above yields the result.
\ \ $\diamond$

\bigskip\noi
{\bf 3.8. Remark: } In [19, Theorem~12.2],
M.~Saito, B.~Sturmfels and N.~Takayama
prove Theorem~3.7 by the method of Gr\"obner deformations.
They also show ([19, Theorem~11.1]) that the lower
bound ${\rm vol}(P) \leq r(\aa)$ holds for arbitrary
$\AA$.

\bigskip

Given $\aa\in \Z^2$, define $\hat\HH(\aa)  := \HH(\aa) /\RR(\aa)$ if
$\RR(\aa) \not= \{0\}$ and 
 $\hat\HH(\aa)  := \HH(\aa) /\C\cdot\chi(\aa)$ if $\RR(\aa) = \{0\}$.

\bigskip
\proclaim Corollary 3.9.   For every $\aa\in\Z^2$,
\par
\item {i)} \ $\dim(\hat\HH(\aa) ) = d-1$.
\smallskip
\item {ii)} \ For every $u\in\N^{m+2}$ the operator
$$D_u\,\colon\, \hat\HH(\aa)  \to \hat\HH(\aa-\AA\cdot u)$$

is an isomorphism.\ \ $\diamond$

\beginsection{\bf References}
 
\bigskip

\item{[1]} S.~Abramov and K.~Kvasenko, 
Fast Algorithms to Search for the Rational Solutions of Linear
Differential Equations with Polynomial Coefficients, {\it in} 
Proceedings of the 1991 International Symposium on
Symbolic and Algebraic Computation (ISSAC'91) (S. M.
Watt, ed.), ACM Press, 1991, 267--270.

\smallskip
\item{[2]} A.~Adolphson, { Hypergeometric functions
and rings generated by monomials}. {\it Duke Math.~J.~} 
{\bf 73} (1994), 269--290.

\smallskip
\item{[3]} A.C.~Avram, E.~Derrick, and D.~Jan\u ci\'c, { On
semi-periods}. 
{\it Nuclear Phys. B} {\bf 471} (1996), 293--308.

\smallskip
\item{[4]} V.~Batyrev, {Variations of the mixed Hodge structure of
affine hypersurfaces in algebraic tori.} {\it Duke Math.~J.~} {\bf 69}
(1993), 349--409.

\smallskip
\item{[5]} V.~Batyrev and D.~van Straten, { Generalized
hypergeometric functions and rational curves on Calabi-Yau complete
intersections in toric varieties}. {\it Commun. Math. Phys.} {\bf 168} (1995),
493--533.

\smallskip
\item{[6]} W.~Bruns and J.~Herzog, {``Cohen--Macaulay Rings,"}
Cambridge Univ. Press, Cambridge, 1993.

\smallskip
\item{[7]} I. M.~Gel'fand, M.~Kapranov, and A.~Zelevinsky,
{Generalized Euler integrals and ${\cal A}$-hypergeometric
functions}. {\it  Advances in Mathematics,.}  {\bf 84} (1990), 255--271.

\smallskip
\item{[8]} I. M.~Gel'fand, M.~Kapranov, and A.~Zelevinsky,
{``Discriminants, Resultants, and Multidimensional Determinants,"}
Birkh\"auser, Boston, 1994.

\smallskip
\item{[9]} I. M.~Gel'fand, A.~Zelevinsky, and  M.~Kapranov,
{Hypergeometric functions and toral manifolds}.
{\it  Functional Analysis and its Applications}
{\bf 23} (1989), 94--106. 

\smallskip
\item{[10]} I. M.~Gel'fand, A.~Zelevinsky, and  M.~Kapranov,
Correction to:{ ``Hypergeometric functions and toral manifolds}." 
{\it  Functional Analysis and its Applications}
{\bf 27} (1994), 295.
 
\smallskip
\item{[11]} I. M.~Gel'fand, A.~Zelevinsky, and  M.~Kapranov,
{$A$-discriminants and Cayley-Koszul complexes}.
{\it  Soviet Math. Doklady }
{\bf 40} (1990), 239--243.
 
\smallskip
\item{[12]} S. Goto, N. Suzuki, and K. Watanabe, { On affine semigroup
rings}. {\it Japan J. Math } {\bf 2} (1976), 1--12.

\smallskip
\item {[13]} S. Hosono, { GKZ systems, Gr\"obner fans and
moduli spaces of Calabi-Yau hypersurfaces}. Preprint. alg-geom/9707003.

\smallskip
\item {[14]} S. Hosono, B.H. Lian, and S.-T.~Yau, { GKZ-Generalized
hypergeometric systems in mirror symmetry of Calabi-Yau 
hypersurfaces}.
{\it Comm.
Math. Phys.} {\bf 182} (1996), 535--577.

\smallskip
\item {[15]} S. Hosono, B.H. Lian, and S.-T.~Yau, { Maximal degeneracy
points of GKZ systems}. {\it Journal of the American Mathematical 
Society} {\bf 10} (1997), 427--443.

\smallskip
\item{[16]} A.~Khovanskii, 
{Newton polyhedra and the Euler-Jacobi formula}, 
{\it Russian Math. Surveys} {\bf 33} (1978), 237--238.

\smallskip
\item{[17]} K.~Mayr, {\"Uber die Aufl\"osung algebraischer
Gleichungssysteme durch hypergeometri\-sche Funktionen}.
{\it Monatshefte f\"ur Mathematik und Physik }
{\bf 45} (1937), 280--313.

\smallskip
\item{[18]}
M.~Saito, B.~Sturmfels, and N.~Takayama, 
{Hypergeometric polynomials and integer programming}.
{\it Compositio Mathematica.} To appear.

\smallskip
\item{[19]}
M.~Saito, B.~Sturmfels, and N.~Takayama, 
{Gr\"obner deformations of hypergeometric
differential equations}.
Preprint.

\smallskip
\item{[20]}{M.~F.~Singer},
{Liouvillian Solutions of $n$-th order
homogeneous linear differential equations\/}. 
{\it American Journal of Mathematics } {\bf 103}  {(1981)},  {661--682}.

\smallskip
\item {[21]} J. Stienstra, { Resonant hypergeometric
systems and mirror symmetry}. Preprint. \hfill \break
 alg-geom/9711002.

\smallskip
\item{[22]} B.~Sturmfels, {Solving algebraic equations in terms
of $\AA$-hypergeometric series}. \hfill \break
{\it  Discrete Mathematics.}  To appear.

\smallskip
\item{[23]}
B.~Sturmfels and N.~Takayama,
{Gr\"obner bases and hypergeometric functions}.  \hfill \break
   {\it Gr\"obner Bases and Applications (Proc. of the Conference
    33 Years of Gr\"obner Bases)},
    B.Buchberger and F.Winkler (eds.),
    Cambridge University Press,
London Mathematical Society Lecture Notes Series, {\bf 251} (1998),
246--258.

\smallskip
\item{[24]} 
N.~Takayama,
{Kan: A system for computation in algebraic analysis}, 
source code available at
{\tt  http://www.math.s.kobe-u.ac.jp/KAN/}  (1991)

\smallskip
\item{[25]} 
N.V. Trung and L.T. Hoa, Affine semigroups and Cohen-Macaulay
rings generated by monomials.  {\it Trans. Amer. Math. Soc.}
{\bf 298} (1986), 145--167.

\bigskip
\noi{\bf Current Addresses:}
\medskip
\noi Eduardo Cattani

\noi Department of Mathematics and Statistics

\noi University of Massachusetts

\noi Amherst, MA 01003

\noi{\tt cattani@math.umass.edu}

\medskip

\noi Carlos D'Andrea

\noi Departmento de Matem\'atica,  F. C. E. y N.

\noi Universidad de Buenos Aires

\noi Ciudad Universitaria -- Pabell\'on I

\noi (1428) Buenos Aires, Argentina.

\noi{\tt cdandrea@dm.uba.ar}

\medskip

\noi Alicia Dickenstein

\noi Departmento de Matem\'atica,  F. C. E. y N.

\noi Universidad de Buenos Aires

\noi Ciudad Universitaria -- Pabell\'on I

\noi (1428) Buenos Aires, Argentina.

\noi{\tt alidick@dm.uba.ar}

\bye